\documentclass[preprint,12pt]{elsarticle}

\newtheorem{thm}{Theorem}
\newtheorem{lem}[thm]{Lemma}
\newtheorem{prop}[thm]{Proposition}
\newtheorem{cor}[thm]{Corollary}
\newdefinition{definition}{Definition}
\newdefinition{rmk}{Remark}
\newdefinition{exm}{Example}
\newdefinition{ass}{Assumption}
\newproof{pf}{Proof}

\usepackage{amsmath,amssymb,amsbsy,amsfonts}
\allowdisplaybreaks

\usepackage{mathtools}
\usepackage{xfrac}
\usepackage{mathrsfs}
\usepackage{enumitem}
\usepackage{dsfont}
\usepackage{bm}

\DeclareMathOperator*{\esslim}{ess\,lim}
\def\fp{ : }

\def\tens#1{\pmb{\mathsf{#1}}}
\def\vec#1{\boldsymbol{#1}}

\def\R{\mathbb{R}}
\def\sym{\mathop{\mathrm{sym}}\nolimits}
\def\Rds{\mathbb{R}^{d \times d}_{\sym}}

\def\tr{\mathop{\mathrm{tr}}\nolimits}
\def\diver{\mathop{\mathrm{div}}\nolimits} %\divergence

\def\wconv{\rightharpoonup}

\def\bf{\vec{f}}

\def\bu{\vec{u}}
\def\bv{\vec{v}}
\def\bw{\vec{w}}

\def\bq{\vec{q}}

\def\BD{\tens{D}}

\def\BG{\tens{G}}

\def\BS{\tens{S}}

\def\Btau{\tens{\tau}}
\def\Bsigma{\tens{\sigma}}
\def\Scr{\bm{\mathcal{S}}}
\def\Dcr{\bm{\mathcal{D}}}
\def\Vndiv{V^n_{\diver}}

\usepackage{color}
\definecolor{morado}{rgb}{0.65, 0.14, 0.78}

\def\texblk#1{\textcolor{black}{#1}}
\journal{Nonlinear Analysis: Real World Applications}

\begin{document}

\begin{frontmatter}



  \title{Weak-strong Uniqueness for Heat Conducting non-Newtonian Incompressible Fluids \tnoteref{t1}}
  \tnotetext[t1]{P. A. Gazca-Orozco's work was supported by the Alexander von Humboldt Stiftung and by the OP RDE project No. CZ.02.2.69/0.0/0.0/16\_027/0008495, International Mobility of Researchers at Charles University. V.\ Patel is supported by the UK Engineering and Physical Sciences Research Council [EP/L015811/1].}

\author[1,2]{Pablo Alexei Gazca-Orozco}
\ead{gazcaorozco@karlin.mff.cuni.cz}
\affiliation[1]{organization={Department of Data Science, FAU Erlangen-Nuernberg},%Department and Organization
            addressline={Cauerstraße 11}, 
            postcode={91058}, 
            city={Erlangen},
            country={Germany}}
\affiliation[2]{organization={Charles University, Faculty of Mathematics and Physics, Mathematical Institute},%Department and Organization
            addressline={{Sokolovska 83}}, 
            postcode={186 75}, 
            city={Prague},
            country={Czech Republic}}
\author[3]{Victoria Patel}
\ead{victoria.patel@maths.ox.ac.uk}
\affiliation[3]{organization={Mathematical Institute, University of Oxford},%Department and Organization
            addressline={Andrew Wiles Building, Woodstock Road}, 
            city={Oxford},
            postcode={OX2 6GG}, 
            country={United Kingdom}}
\begin{abstract}
  In this work, we introduce a notion of dissipative weak solution for a system describing the evolution of a heat-conducting incompressible non-Newtonian fluid. This concept of solution is based on the balance of entropy instead of the balance of energy and has the advantage that it admits a weak-strong uniqueness principle, justifying the proposed formulation. We provide a proof of existence of solutions based on finite element approximations, thus obtaining the first convergence result of a numerical scheme for the full evolutionary system including temperature dependent coefficients and viscous dissipation terms. Then we  proceed to prove the weak-strong uniqueness property of the system by means of a relative energy inequality. 
\end{abstract}




\end{frontmatter}


\section{Introduction and problem formulation}
%\subsection*{Statement of the problem}
For \( d\in \{ 2,3\}\), let \( \Omega\subset\mathbb{R}^d \) be a bounded, Lipschitz domain. For a given final time horizon \( T
\in (0, \infty) \), we define the space-time domain \( Q = (0, T) \times\Omega\). % For any \( \tau \in (0, T) \), we denote the parabolic cylinder \( Q_\tau = (0, \tau)\times \Omega\).
Given a body force \( \bm{f}\colon Q \rightarrow\mathbb{R}^d\), an initial velocity field \( \bu_0 \colon \Omega \rightarrow\mathbb{R}^d\) and an initial internal energy \( e_0 \colon\Omega\rightarrow\mathbb{R}\), we consider the problem of finding a divergence free velocity field \( \bu\colon \overline{Q}\rightarrow\mathbb{R}^d\), a positive internal energy \( e\colon\overline{Q}
\rightarrow\mathbb{R}\), a pressure \( p \colon Q \rightarrow\mathbb{R}\), \texblk{a heat flux $\bm{q}\colon Q\to \R^d$}, and a traceless stress tensor field \( \BS \colon Q\rightarrow\Rds\) such that we have the balance laws

\begin{subequations}\label{eq:PDE}
\begin{alignat}{2}
%\begin{aligned}
%\label{eq:PDE}
%<<<<<<< HEAD
%\partial_t e + \diver(e\bu\, +\, &\bq) = \BS \fp \BD(\bu)  &&\quad \text{ on } Q,\label{eq:PDE_energy}\\
%\partial_t \bu - \diver(\BS - \bu &\otimes  \bu)
%+\nabla p = \bf \qquad \quad &&\quad \text{ on } Q,\label{eq:PDE_momentum}\\
% \diver \bu &= 0  \qquad \quad &&\quad \text{ on }{Q},\label{eq:PDE_mass}
%=======
\partial_t \bu - \diver(\BS - \bu &\otimes  \bu)
+\nabla p = \bf \qquad \quad &&\quad \text{ on } Q,\label{eq:PDE_momentum}\\
 \diver  &\bu = 0  \qquad \quad &&\quad \text{ on }{Q},\label{eq:PDE_mass}\\
\partial_t e + \diver(e\bu\, +&\,\, \bq) = \BS \fp \BD\bu  &&\quad \text{ on } Q,\label{eq:PDE_energy}
%\end{aligned}
\end{alignat}
subject to the  initial and boundary conditions
\begin{alignat}{2}
\begin{aligned}\label{eq:PDE_bcs}
\bu(0,\cdot)&=\bu_0(\cdot)  \qquad &&\text{ in } \Omega,\\
e(0,\cdot)&=e_0(\cdot)  \qquad &&\text{ in } \Omega, \\
\bu &= \bm{0}  \qquad \quad &&\text{ on } (0,T) \times \partial \Omega, \\
\bq \cdot \bm{n} &= 0  \qquad \quad &&\text{ on } (0,T) \times \partial \Omega.
\end{aligned}
\end{alignat}
%\end{subequations}
%
We  assume that the internal energy is related to the temperature by  $e= c_v\theta$ where $c_v>0$ is assumed to be a constant.
The system  is closed by relating the heat flux $\bm{q}\colon Q \to \mathbb{R}^d$ and the stress tensor $\BS$ to the temperature gradient $\nabla \theta$ and the symmetric velocity gradient $\BD\bu := \frac{1}{2}(\nabla\bu + \nabla\bu^\top)$, respectively, through constitutive relations of the form
%\begin{subequations}\label{eq:constitutive_rel}
\begin{alignat}{2}
  \bq &= -\kappa(e)\nabla e = -\tilde{\kappa}(\theta) \nabla \theta \qquad &&\text{ a.e. in }Q,\\
%    (\BD(\bu),\BS,e)&\in \mathcal{A}(\cdot) &&\text{ a.e. in }Q.
\BS &= \Scr(\BD\bu,\theta) && \text{ a.e. in }Q,
\end{alignat}
\end{subequations}
where $\Scr\colon \Rds\times \R \to \Rds$ and $\tilde{\kappa}\colon \R \to \R$ are given continuous functions. The precise assumptions will be introduced in Section \ref{sec:exist}.

One of the main challenges in the analysis of system \eqref{eq:PDE} arises from the presence of the viscous dissipation term $\BS \fp \BD\bu$ in the balance of internal energy \eqref{eq:PDE_energy}. The difficulty stems from the fact that this term belongs a priori only to $L^1(Q)$, which makes the application of compactness arguments problematic. For this reason, most of the early works that tackled the question of existence of solutions for systems describing incompressible heat-conducting fluids either neglected viscous dissipation \cite{Kagei1993,Malek1994,Hishida1992} or employed a weaker notion of weak solution such as a distributional solution or weak solution with a defect measure \cite{Diaz1998,Necas2001,Naumann2006}. The works \cite{Clopeau1997,Consiglieri2000} instead employed a setting in which the velocity $\bu$ is an admissible test function in the weak formulation of the balance of momentum \eqref{eq:PDE_momentum}, which simplifies some of the arguments, but excludes the Navier--Stokes model in three dimensions.

A breakthrough came with the work of \texblk{Bul\'{i}\v{c}ek, Feireisl and M\'{a}lek} \cite{Bulicek2009} (see also \cite{Feireisl2006}),  where it was observed that, even though it contains additional couplings, a formulation that employs the following equation for the total energy $E:= \tfrac{1}{2}|\bu|^2 + e$ instead of the balance of internal energy \eqref{eq:PDE_energy} is more amenable to weak convergence arguments:
\begin{equation}\label{eq:PDE_total_energy}
    \partial_t E + \diver((E + p)\bu - \BS\bu) - \diver(\kappa(e)\nabla e) = \bm{f}\cdot\bu. %+ g.
\end{equation}
In particular, the existence of bona fide weak solutions for the system with Newtonian rheology and temperature dependent coefficients was established for large data. We note that the two formulations \eqref{eq:PDE_total_energy} and \eqref{eq:PDE_energy} are equivalent when the solutions are smooth or, in the weak setting, whenever it is allowed to test the momentum balance with the velocity $\bu$. This idea was further applied in \cite{Bulicek2009a,Maringova2018} to models with shear-rate and pressure dependent viscosities and implicit models with activation parameters, respectively. A drawback of the formulation involving \eqref{eq:PDE_total_energy} is that one needs an integrable pressure, which precludes the use of the popular no-slip boundary condition $\bu|_{\partial\Omega}=\bm{0}$ for the velocity (e.g.\ Navier's slip boundary condition was employed in \cite{Bulicek2009,Bulicek2009a}). Furthermore, the approach in \cite{Bulicek2009,Maringova2018} makes use of regularity properties of the Neumann--Laplace problem when obtaining {\it a priori} estimates for the pressure which requires more than mere Lipschitz regularity of the domain $\Omega$.
%It allows to passage to the limit more easily when considering the limit of approximations of the system \eqref{eq:PDE}. Moreover, \eqref{eq:PDE_energy} is recovered as an inequality. For smooth solutions,  the formulations are in fact equivalent.  The disadvantage of this formulation is that one needs an integrable pressure, which requires a slip boundary condition and higher regularity of the domain. Motivated by a desire to have a problem formulation that allows for Lipschitz domains, we avoid the use of (\ref{eq:PDE_total_energy}) and the total energy. Instead, we look for an equivalent form of (\ref{eq:PDE_energy}) that does not involve a pressure term.

In this work we follow an alternative approach and introduce a notion of \emph{dissipative weak solution} to the system with no-slip boundary conditions for the velocity on general Lipschitz domains (hence including polyhedral/polygonal domains usually employed in numerical approximations). The formulation here is inspired from the works \cite{Feireisl2012,Feireisl2012A,Brezina2018}, which dealt with the compressible Navier--Stokes--Fourier system, and is based on the balance of entropy rather than the balances of energy \eqref{eq:PDE_energy} or \eqref{eq:PDE_total_energy}. We define the entropy \( S = c_v \log\theta\) and consider the entropy balance
\begin{equation}\label{eq:entropy_smooth}
    \partial_t S + \diver (S\bu) + \diver \left(\frac{\bq}{\theta}\right)
    \geq \frac{1}{\theta}\left( \BS\colon\BD\bu - \frac{\bq}{\theta}\cdot \nabla \theta \right),
\end{equation}
supplemented with the total energy balance
\begin{equation}\label{eq:total_energy_smooth}
  \frac{\mathrm{d}}{\mathrm{d}t} \int_\Omega\left[ \frac{1}{2}|\bu|^2 + e \right]
  \leq \int_\Omega \bm{f}\cdot\bu .
\end{equation}
We note that for smooth solutions, the formulations with \eqref{eq:PDE_energy} or \eqref{eq:PDE_total_energy} are equivalent to the formulation with \eqref{eq:entropy_smooth}, with an equality sign ``\(=\)'' replacing the inequality sign ``\(\geq \)''. The balance \eqref{eq:total_energy_smooth} is in turn also satisfied as an equality. The following definition states precisely the concept of solution that we consider in this work. The notation employed here is properly defined in the next section.

\begin{definition}\label{def:dissipative_weak_sols}
  Let $r>\tfrac{2d}{d+2}$ and suppose that we are given a function $\bm{f}\in L^{r'}(0,T;W^{-1,r'}(\Omega)^d)$.  We say that a triple \((\BS, \bu, \theta) \) is a dissipative weak solution of (\ref{eq:PDE}) if
\begin{align*}
  \BS &\in L_{\sym,\tr}^{r^\prime}(Q)^{d\times d},
\\
  \bu &\in L^r(0, T ; W_{0,\diver}^{1,r}(\Omega)^d) \cap L^\infty(0, T ; L_{\diver}^2(\Omega)^d),
  \\
%  \partial_t \bu &\in L^{\check{r}}(0, T; (W^{-1,\check{r}'}_{0,\diver}(\Omega)^d)^*),
%  \\
\theta&\in L^{q_1}(0, T; W^{1,q_1}(\Omega)) \cap L^{q_2}(Q),
\\
\log\theta&\in L^2(0, T; W^{1,2}(\Omega)) \cap L^\infty(0, T; L^{q_3}(\Omega)),
\end{align*}
 \texblk{for any} \( q_1 \in [1,\frac{5}{4}) \), \( q_2 \in [1,\frac{5}{3}) \), \( q_3 \in [1,\infty) \) and the following relations are satisfied:
\begin{subequations}\label{eq:weakPDE}
\begin{itemize}[leftmargin=2.3ex]
  \item the constitutive relation holds pointwise almost everywhere,
\begin{equation}\label{eq:weakPDE_CR}
  \BS = \Scr (\BD\bu,\theta)\quad \text{a.e. in }Q,
\end{equation}
where \( \Scr: \mathbb{R}^{d\times d}_{\sym}\times\mathbb{R}\rightarrow\mathbb{R}^{d\times d}_{\sym}\) is a given continuous function;
\item the balance of momentum holds in the usual weak sense, 
%  \begin{equation}\label{eq:weakPDE_momentum}
%\int_{Q_\tau} \bu \cdot \partial_t \bv\,\mathrm{d}x\,\mathrm{d}t - \Big[ \int_\Omega\bu(t) \cdot \bv(t) \,\mathrm{d}x
%\Big]_{t = 0}^{t = \tau}
%%\\
%= \int_{Q_\tau} (\BS - \bu \otimes \bu ) : \BD \bv \,\mathrm{d}x\,\mathrm{d}t - \int_0^\tau\langle \bf, \bv \rangle\,\mathrm{d}t,
%\end{equation}
  \begin{equation}\label{eq:weakPDE_momentum}
\int_{Q_\tau} \bu \cdot \partial_t \bv - \Big[ \int_\Omega\bu(t) \cdot \bv(t) 
\Big]_{t = 0}^{t = \tau}
%\\
= \int_{Q_\tau} (\BS - \bu \otimes \bu ) \fp \BD \bv 
- \int_0^\tau\langle \bf, \bv \rangle,
\end{equation}
for every $\bv \in C^\infty_0([0,T);C^\infty_{0,\diver}(\Omega)^d)$ and a.e. \( \tau \in (0, T) \);
\item for the entropy $S:= c_v \log \theta$, the entropy inequality holds weakly in the sense that
%  \begin{equation}\label{eq:weakPDE_entropy}
%\begin{aligned}
%&-\int_{Q_\tau} S\partial_t \psi \,\mathrm{d}x\,\mathrm{d}t + \Big[ \int_\Omega \psi(t)c_v \log\theta(t) \,\mathrm{d}x\Big]_{t = 0}^{t = \tau} - \int_{Q_\tau} S\bu \cdot \nabla \psi \,\mathrm{d}x\,\mathrm{d}t
%\\
%&
%+ \int_{Q_\tau} \frac{\tilde{\kappa}(\theta) \nabla\theta}{\theta} \cdot \nabla \psi \,\mathrm{d}x\,\mathrm{d}t
%\geq \int_{Q_\tau} \frac{\tilde{\kappa}(\theta) |\nabla \theta|^2}{\theta^2 }\psi \,\mathrm{d}x\,\mathrm{d}t + \int_{Q_\tau} \frac{\BS : \BD \bu }{\theta}\psi \,\mathrm{d}x\,\mathrm{d}t,
%\end{aligned}
%\end{equation}
\begin{equation}\label{eq:weakPDE_entropy}
\begin{aligned}
&-\int_{Q_\tau} S\partial_t \psi 
+ \Big[ \int_\Omega \psi(t) S(t) \Big]_{t = 0}^{t = \tau} 
- \int_{Q_\tau} S\bu \cdot \nabla \psi 
\\
&
+ \int_{Q_\tau} \frac{\tilde{\kappa}(\theta) \nabla\theta}{\theta} \cdot \nabla \psi 
\geq \int_{Q_\tau} \frac{\tilde{\kappa}(\theta) |\nabla \theta|^2}{\theta^2 }\psi  
+ \int_{Q_\tau} \frac{\BS \fp \BD \bu }{\theta}\psi ,
\end{aligned}
\end{equation}
for any \( \psi \in C^\infty_0([0, T) ; C^\infty(\overline{\Omega})) \) such that \( \psi \geq 0 \) and a.e. \( \tau \in (0, T) \);
%{\color{blue}
%\item the incompressibility constraint holds in the usual weak sense,
%  \begin{equation}\label{eq:weakPDE_mass}
%\int_{Q} q \diver \bu = 0, 
%  \end{equation}
%  for every $q\in C_0^\infty(Q)$;
%}
\item the total energy inequality holds, 
\end{itemize}
%\begin{equation}\label{eq:weakPDE_energy}
%\Big[ \int_\Omega   \frac{|\bu(t,\cdot) |^2}{2} + \theta(t,\cdot)   \,\mathrm{d}x\Big]_{t = 0}^{t = \tau} \leq \int_0^\tau \langle\bf, \bu\rangle \,\mathrm{d}t,
%\end{equation}
\begin{equation}\label{eq:weakPDE_energy}
\Big[ \int_\Omega   \frac{|\bu(t,\cdot) |^2}{2} + \theta(t,\cdot)   \Big]_{t = 0}^{t = \tau} 
\leq \int_0^\tau \langle\bf, \bu\rangle ,
\end{equation}
for a.e. $\tau\in (0,T)$.
\end{subequations}
\end{definition}

%\AG{I'm  curious about something Josef asked: if one assumes additional regularity, does our dissipative weak solution coincide with the weak solution from their works? I would imagine that if $r>\frac{3d}{d+2}$ and the pressure is a Lebesgue function, then this is probably true...}
%\VP{I'm not entirely sure. I can definitely believe it but I think it might take some work to prove rigorously?}

The parameter \( r\) is determined by the coercivity property satisfied by the constitutive relationship (see Assumption \ref{as:const_rel} below); for instance, for the Navier--Stokes model one has $r=2$. \texblk{The restriction $r>\frac{2d}{d+2}$ is a natural one, since it guarantees that the convective term in the momentum equation can be handled as a compact perturbation, thanks to the compact embedding of $W^{1,r}(\Omega)^d$ into $L^2(\Omega)^d$}.

The advantage of such a formulation involving the entropy is that the corresponding solutions satisfy a \emph{weak-strong uniqueness principle}, i.e.,  the dissipative weak solution will be equal to the strong solution emanating from the same initial data for as long as the latter exists. The fact that a weak-strong uniqueness result holds is an indicator that the notion of weak solution under consideration is a sensible extension of the classical one, and thus the result is of interest on its own right. However, weak-strong uniqueness results can also be useful in the analysis of singular limits and stability of stationary states \cite{Saint-Raymond2009,Brezina2018} and have been obtained in different contexts \cite{Feireisl2012,Demoulini2012,Lattanzio2013,Abbatiello2019a}.

In the present work, we first prove in Section \ref{sec:exist} the existence of solutions to the system \eqref{eq:PDE} in the sense of Definition \ref{def:dissipative_weak_sols}. The existence proof employs similar ideas to the ones presented in \cite{Bulicek2009}, with a couple of important differences. In \cite{Bulicek2009} an abstract Galerkin approach is first applied to the system using a quasi-compressible approximation $\varepsilon\Delta p_\varepsilon = \diver \bu_\varepsilon$. Additionally, the convective term is handled by using a divergence-free mollifier approximate $\bu_\varepsilon$ (constructed with the help of a Helmholtz decomposition). In contrast, in this work we construct the approximations by means of a numerical scheme based on the finite element method using the usual divergence-free constraint $\diver \bu = 0$ and standard LBB \texblk{(Ladyzhenskaya--Babu\v{s}ka--Brezzi)} stable finite element spaces. Furthermore, the convective terms can be handled using the typical skew-symmetric form employed in numerical analysis, thus avoiding the use of a Helmholtz decomposition and mollifiers, which would complicate the implementation of the numerical scheme. Since the formulation considered here does not involve the balance \eqref{eq:PDE_total_energy}, our result can be obtained by assuming $r>\frac{2d}{d+2}$, which is the natural assumption required to handle the convective term $\diver (\bu\otimes \bu)$ as a compact perturbation, and is less restrictive than the condition $r>\frac{3d}{d+2}$, which was needed in the works \cite{Bulicek2009,Bulicek2009a}. We should also mention that while Definition \ref{def:dissipative_weak_sols} considers an explicit constitutive relation \eqref{eq:weakPDE_CR}, the approach employed here is well suited to handle models with implicit constitutive relations (see Remark \ref{rem:implicitCR}).

In Section \ref{sec:WSU}, we proceed to prove that a corresponding weak-strong uniqueness principle applies to our notion of dissipative weak solution using the method of relative entropies. This could be considered as the incompressible non-Newtonian counterpart of the results from \cite{Feireisl2012,Feireisl2012A,Brezina2018}. For the incompressible Navier--Stokes model one has the classical results of Prodi \cite{Prodi1959} and Serrin \cite{Serrin1962} (see also \cite{Wiedemann2018} for a more recent survey) and, more recently, a weak-strong uniqueness result for the incompressible non-Newtonian system with an implicit constitutive relation was obtained in \cite{Abbatiello2019}. The present work can be considered an extension of \cite{Abbatiello2019} to the non-isothermal setting with temperature-dependent coefficients.

Finally, we  highlight the fact that since the dissipative weak solutions were constructed by means of a numerical scheme, as a consequence we   obtain here the first convergence result of finite element approximations to a solution of a system describing a heat-conducting non-Newtonian incompressible fluid with no-slip boundary conditions for the velocity, using a model that does not neglect viscous dissipation. Similar ideas can be found in \cite{Feireisl2019,Li2021}, where convergence of certain finite volume schemes was established for some compressible fluid models.

\section{Finite element solutions generate dissipative weak solutions}\label{sec:exist}
Throughout this work, we  employ  standard notation for Lebesgue,  Sobolev and Bochner  spaces (e.g.\ $(W^{k,p}(\Omega), \|\cdot\|_{W^{k,r}(\Omega)})$ and $(L^q(0,T;W^{k,p}(\Omega), \|\cdot\|_{L^q(0,T;W^{k,p}(\Omega))}))$). The space $W^{k,p}_0(\Omega)$ is defined as the closure of the space of smooth and compactly supported functions $C_0^\infty(\Omega)$ with respect to the $\|\cdot\|_{W^{k,p}(\Omega)}$ norm. \texblk{We denote their divergence-free subspaces as $W^{k,p}_{0,\diver}(\Omega)^d := \{\bv\in W^{1,r}_{0}(\Omega)^d \,:\, \diver \bv=0\}$ and $C^\infty_{0,\diver}(\Omega)^d := \{\bv\in C_0^\infty(\Omega)^d \, :\, \diver \bv=0\}$}. 
The dual space of $W_0^{1,r}(\Omega)$ is denoted by $W^{-1,r'}(\Omega)$, where $r'$ is the H\"{o}lder conjugate of $r\in (1,\infty)$. Finally, the space $L^q(Q)^{d\times d}_{\sym,\tr}$ will denote the subspace of matrix-valued functions in $L^q(Q)^{d\times d}$ that are symmetric and traceless.

In this section, we prove the existence of dissipative weak solutions to the system (\ref{eq:PDE}) in the sense of Definition \ref{def:dissipative_weak_sols}. In order to proceed,  we need appropriate monotonicity and coercivity assumptions on the constitutive relation.

\begin{ass}\label{as:const_rel}
  The function $\Scr\colon \Rds\times \R \to \Rds$ defining the constitutive relation \eqref{eq:weakPDE_CR} is continuous and satisfies the following further properties.
\begin{itemize}[leftmargin = 0.8cm]
	\item {\it (Monotonicity)} For every fixed \( s\in \R\) and for every $\Btau_1,\Btau_2 \in \Rds$,
		\begin{equation}\label{eq:monotonicity}
			(\Scr(\Btau_1,s)-\Scr (\Btau_2,s))\fp (\Btau_1 - \Btau_2)\geq 0.
		\end{equation}
	\item {\it (Coercivity)} There exist a non-negative function $g\in L^{1}(Q)$ and a constant $c>0$ such that
    \begin{equation}\label{eq:coercivity}
			\Scr(\Btau,s)\fp \Btau \geq -g + c(|\Scr(\Btau,s)|^{r'} + |\Btau|^r)\quad \text{ for all }\Btau\in \Rds\text{ and }s\in \R.
		\end{equation}
	\item {\it (Growth)} There exists a constant $c>0$ such that
    \begin{equation}\label{eq:growth}
			|\Scr(\Btau,s)| \leq c(|\Btau|^{r-1} + 1) \quad \text{ for all }\Btau\in \Rds\text{ and } s\in \R.
		\end{equation}
	\item {\it (Compatibility)} For every fixed $s\in\R$, we have    $\tr(\Scr(\Btau,s)) = 0$ if and only if $\tr(\Btau) = 0$, for every $\Btau\in \Rds$.
\end{itemize}
Regarding the heat flux, we assume that the heat conductivity $\tilde{\kappa}\colon \R\to \R$ is a continuous function such that $0<c_1 \leq \tilde{\kappa}(s) \leq c_2$, for any $s\in \R$, where $c_1$ and $c_2$ are two positive constants.
\end{ass}

Under these assumptions, the focus of this section is to prove the existence of weak dissipative solutions of (\ref{eq:PDE}).
The proof is based on a \( 3\)-level approximation scheme, motivated by finite element techniques used in numerical analysis. The approximation indices are
\begin{itemize}
\item \(m\) for the time discretisation and the Galerkin discretisation for \( \theta\),
\item \( n \) for the Galerkin discretisation for \( \bu \), and
\item \( k \) for a penalty term. %{\color{blue} This is also needed as a cut-off for the constitutive relation I think.}
\end{itemize}
The presence of a penalty term allows us to test in the momentum balance against the velocity itself, despite the presence of the quadratic term \( \bu \otimes \bu \). We look for a solution of the \(3\)-level approximation scheme. Then we take the limit in \(m\), followed by \( n \) and finally \(k\) (see also Remark \ref{rem:3lev_approx}).
For notational simplicity, we denote \( \alpha = (m,n,k)\) and \( \beta = (n,k) \).

Let $\{\mathcal{T}_n\}_{n\in\mathbb{N}}$ be a family of shape-regular triangulations of $\Omega$, for which the mesh size $h_n := \max_{K\in \mathcal{T}_n}\mathrm{diam}(K)$ vanishes as $n\to\infty$. We define the following conforming finite element spaces for the temperature, velocity and pressure, respectively:
\begin{align*}
U^n &:= \{w \in W^{1,\infty}(\Omega) \, :\, w|_K \in \mathbb{P}_{\mathbb{U}}(K), \, K\in \mathcal{T}_n\},\\
  V^n &:= \{\bv \in W_0^{1,\infty}(\Omega)^d \, :\, \bv|_K \in \mathbb{P}_{\mathbb{V}}(K)^d, \, K\in \mathcal{T}_n\},\\
  M^n &:= \{q \in L^\infty(\Omega) \, :\, q|_K \in \mathbb{P}_{\mathbb{M}}(K), \, K\in \mathcal{T}_n\},
\end{align*}
where $\mathbb{P}_{\mathbb{U}}(K)$, $\mathbb{P}_{\mathbb{V}}(K)$, and $\mathbb{P}_{\mathbb{M}}(K)$ denote spaces of polynomials on the element $K\in \mathcal{T}_n$. These must be chosen in such a way that certain stability properties are satisfied (see Assumption \ref{as:FEM} below). We also introduce the following useful subspace of discretely divergence free functions of $V^n$:
\begin{equation*}
\Vndiv := \left\{\bv \in V^n \,  :\, \int_\Omega q \,\diver \bv =0\quad \text{for all }q\in M^n \right\}.
\end{equation*}
Although the formulation \eqref{eq:weakPDE} does not include the pressure, in practice the incompressibility constraint is enforced at the discrete level by means of a Lagrange multiplier, which could be interpreted as a discrete pressure. For this reason, we also introduce appropriate assumptions on the pressure space in what follows.

\begin{ass}\label{as:FEM}
  The finite element spaces $U^n$, $V^n$, and $M^n$ satisfy the following properties.
\begin{itemize}[leftmargin = 0.8cm]
  \item {\it (Approximability)} For an arbitrary $s\in [1,\infty)$, one has that
    \begin{align*}
         \inf_{\overline{w}\in U^n}\|w - \overline{w}\|_{W^{1,s}(\Omega)} &\to 0 &\quad& \textrm{as }n\to\infty \quad \forall\, w\in W^{1,s}(\Omega),\\
      \inf_{\overline{\bv}\in V^n}\|\bv - \overline{\bv}\|_{W^{1,s}(\Omega)} &\to 0 &\quad& \textrm{as }n\to\infty \quad \forall\, \bv\in W^{1,s}(\Omega)^d,\\
      \inf_{\overline{q}\in M^n}\|q - \overline{q}\|_{L^{s}(\Omega)} &\to 0 &\quad &\textrm{as }n\to\infty \quad \forall\, q\in L^{s}(\Omega).
    \end{align*}
  \item {\it (Fortin Projector $\Pi^n_V$)} For every $n\in \mathbb{N}$, there exists  $\Pi^n_V\colon W^{1,1}_0(\Omega)^d \to V^n$, a linear projector,  that satisfies the usual stability and divergence preservation properties. That is,  for any $\bv\in W^{1,1}_0(\Omega)^d$, we have
    \begin{align*}
    \int_\Omega q\,\diver\bv &= \int_\Omega q \,\diver(\Pi^n_V\bv) \quad \forall\, q\in M^n,
\end{align*}
and
\begin{align*}
\|\Pi^n_V\bv\|_{W^{1,s}(\Omega)} &\leq c \|\bv\|_{W^{1,s}(\Omega)},
    \end{align*}
    where $s\in [1,\infty)$ is arbitrary and $c>0$ is a constant that is independent of $n$.
  \item {\it (Projectors $\Pi^n_U$, $\Pi^n_M$)} For every $n\in \mathbb{N}$, we assume that there exist   $\Pi^n_U\colon W^{1,1}(\Omega)\to U^n$ and $\Pi^n_M\colon L^1(\Omega)\to M^n$, linear projectors, such that
    \begin{align*}
    \|\Pi^n_U w\|_{W^{1,s}(\Omega)} &\leq c \|w\|_{W^{1,s}(\Omega)} &\quad &\forall\, w\in W^{1,s}(\Omega),\\
      \|\Pi^n_M q\|_{L^{s}(\Omega)} &\leq c \|q\|_{L^{s}(\Omega)} &\quad &\forall\, q\in L^{s}(\Omega),
    \end{align*}
    where $s\in [1,\infty)$ is arbitrary and $c>0$ is a constant, independent of $n$.
\end{itemize}
\end{ass}

{\color{black}
A direct consequence of the approximability and stability properties in Assumption \ref{as:FEM} is that, for any $s\in [1,\infty)$,
\begin{align*}
\|\bv - \Pi^n_V\bv\|_{W^{1,s}(\Omega)} \to 0 \quad &\text{as} \quad n\to \infty,\notag \\
\|q - \Pi^n_Mq\|_{L^s(\Omega)} \to 0 \quad &\text{as} \quad n\to\infty,\\
\|w - \Pi^n_Uw\|_{W^{1,s}(\Omega)} \to 0 \quad &\text{as} \quad n\to\infty, \notag
\end{align*}
which is useful when passing to the limit in the numerical scheme. There are several known examples of finite element spaces satisfying  Assumption \ref{as:FEM} (see e.g.\ \cite{Boffi2013,Belenki:2012,Zhang2005,Scott1990}).
}

In the discretisation scheme, we   employ the skew-symmetric form of the convective term. More precisely, the trilinear forms meant to represent the convective terms in the momentum and temperature equations are defined, respectively, as
\begin{equation*}
\renewcommand{\arraystretch}{1.5}
\mathcal{B}(\bm{u},\bm{v},\bm{w}) := \left\{
\begin{array}{cc}
  - \displaystyle\int_\Omega\bm{u}\otimes\bm{v}\fp\nabla \bm{w}  & \textrm{ if } V^n_{\diver}
  \subset W^{1,1}_{0,\textrm{div}}(\Omega)^d,\\
  \displaystyle\frac{1}{2}\int_\Omega \left[ \bm{u}\otimes\bm{w}\fp\nabla\bm{v}-\bm{u}\otimes\bm{v}\fp\nabla \bm{w}\right] & \textrm{ otherwise},\\
\end{array}
\right.
\end{equation*}
and
\begin{equation*}
\renewcommand{\arraystretch}{1.5}
\mathcal{C}(\bm{u},\theta,\eta) := \left\{
\begin{array}{cc}
  - \displaystyle\int_\Omega\bm{u}\theta \cdot\nabla \eta  & \textrm{ if } V^n_{\diver} \subset W^{1,1}_{0,\textrm{div}}(\Omega)^d,\\
  \displaystyle\frac{1}{2}\int_\Omega \left[ \bm{u}\eta\cdot\nabla\theta - \bm{u}\theta \cdot \nabla \eta\right] & \textrm{ otherwise}.\\
\end{array}
\right.
\end{equation*}
The advantage of this choice is that we recover the skew-symmetry property that is valid for the original convective term at the continuous level. Indeed, we have  $\mathcal{B}[\bu,\bv,\bv]=0$ and $\mathcal{C}[\bv,\eta,\eta]=0$ for any $\bu$, $\bv\in \Vndiv$ and $\eta\in U^m$, regardless of whether the discretely divergence free velocities are also pointwise divergence free or not. This  is crucial to obtaining {\it a priori} estimates on the sequence of approximate solutions. 

For the time discretisation, we take a sequence of time steps $\{\tau_m\}_{m\in\mathbb{N}}$ such that $T/\tau_m\in \mathbb{N}$ and $\tau_m\to 0$ as $m\to \infty$. For each time step \( \tau_m\), we work on the equidistant grid \( \{t_j\}_{j=0}^{T/\tau_m} \) where we define \( t_j := t^m_j := j\tau_m\) for \( 0 \leq j \leq T/\tau_m\).

Given a sequence of functions $\{v_j\}_{j=0}^{T/\tau_m} $ belonging to some Banach space $X$, we define the piecewise constant interpolant $\overline{v}\in L^\infty(0,T;X)$ by
\begin{subequations}
\begin{equation}\label{def:constant_interpolant}
  \overline{v}(t) := v_j\qquad
  \text{for } t\in (t_{j-1},t_j] ,\,
  j\in \{1,\ldots, T/\tau_m\}.
\end{equation}
The piecewise linear interpolant $\tilde{v}\in C([0,T];X)$ is defined by
\begin{equation}\label{def:linear_interpolant}
  \tilde{v}(t) := \frac{t-t_{j-1}}{\tau_m}v_j + \frac{t_j - t}{\tau_m}v_{j-1}  \qquad
  \text{for }t\in [t_{j-1},t_j],\,
  j\in \{1,\ldots, T/\tau_m\}.
\end{equation}
Additionally, we define the time averages of a given function $g\in L^p(0,T;X)$ by
\begin{equation}\label{def:time_average}
  g_j(\cdot) := \frac{1}{\tau_m}\int_{t_{j-1}}^{t_j} g(t,\cdot) \,\mathrm{d}t.
\end{equation}
\end{subequations}
It is possible to then prove that the piecewise constant interpolant $\overline{g}^m$ defined by \eqref{def:constant_interpolant} for the sequence \((g_j)_{j=0}^{T/\tau_m}\) satisfies $\|\overline{g}^m\|_{L^p(0,T;X)}\leq \|g\|_{L^p(0,T;X)}$ and $\overline{g}^m \to g$ strongly in $L^p(0,T;X)$ as $m\to \infty$ \cite{Roubicek2013}.

The formulation of the discrete problem with parameters \(\alpha = (m,n ,k ) \) and corresponding existence result is as follows.

\begin{prop}\label{exist:prop1}
Let $r>\frac{2d}{d+2}$ and let \( \alpha = (m,n,k) \) be fixed approximation parameters. Suppose that the data \( \bf \in L^{r^\prime}(0, T; W^{-1,r^\prime}(\Omega)^d) \), \( \bu_0 \in L^2_{\diver}(\Omega)^d\) and \( \theta_0 \in L^1(\Omega)\) are given, such that \( \theta_0 \geq c_*> 0 \) for a constant \( c_*\). We define the initialisations
\begin{align*}
\bu^\alpha_0 = P^n_V \bu_0, \quad \theta^\alpha_0 = P^m_U\theta^n_0,
\end{align*}
where \( P^n_V\), \( P^m_U \) are the \( L^2\)-projection operators onto \( V_{\diver}^n\), \( U^m \), respectively, and \(\theta^n_0\) is  defined as follows. We extend \( \theta_0 \) by \( c_* \) outside of \( \Omega\) and    define \( \theta_0^n = \rho_{\frac{1}{n}}* \theta_0\),  where \(\rho_{\frac{1}{n}}\) is a mollification kernel of radius \( \frac{1}{n}\). Let   $\bm{f}^m_j\in W^{-1,r'}(\Omega)^d$ be the sequence of time averages associated to $\bm{f}$. Define $r^\circledast$ to be an arbitrary but fixed number that is greater than $\max\{2r',5\}$.

For every \( j \in \{ 1,\dots, m\} \), defining solutions recursively, there exist \( \bu^\alpha_j \in \Vndiv \) and \( \theta^\alpha_j \in U^m \) such that
%\begin{equation}\label{exist:equ9}
%\int_\Omega \delta\bu^\alpha_j \cdot \bv \,\mathrm{d}x+ \int_\Omega \Scr ( \BD\bu^\alpha_j , \theta^\alpha_j) \fp \BD \bv  \,\mathrm{d}x
%%- (\bu^\alpha_j \otimes \bu^\alpha_j) \Phi_k(|\bu^\alpha_j|) : \BD \bv \,\mathrm{d}x
%%+ \frac{1}{k}\int_\Omega |\bu^\alpha_j|^{2r'-2}\bu^\alpha_j \cdot \bv
%+ \frac{1}{k}\int_\Omega |\bu^\alpha_j|^{r^\circledast-2}\bu^\alpha_j \cdot \bv\,\mathrm{d}x
%+ \mathcal{B}[\bu^\alpha_j, \bu^\alpha_j, \bv]
%= \langle \bm{f}_j^m, \bv\rangle,
%\end{equation}
\begin{equation}\label{exist:equ9}
\int_\Omega \delta\bu^\alpha_j \cdot \bv 
+ \int_\Omega \Scr ( \BD\bu^\alpha_j , \theta^\alpha_j) \fp \BD \bv  
%- (\bu^\alpha_j \otimes \bu^\alpha_j) \Phi_k(|\bu^\alpha_j|) : \BD \bv \,\mathrm{d}x
%+ \frac{1}{k}\int_\Omega |\bu^\alpha_j|^{2r'-2}\bu^\alpha_j \cdot \bv
+ \frac{1}{k}\int_\Omega |\bu^\alpha_j|^{r^\circledast-2}\bu^\alpha_j \cdot \bv
+ \mathcal{B}[\bu^\alpha_j, \bu^\alpha_j, \bv]
= \langle \bm{f}_j^m, \bv\rangle,
\end{equation}
and
%\begin{equation}\label{exist:equ10}
%\int_\Omega\delta\theta^\alpha_j \psi\,\mathrm{d}x
%+ \int_\Omega \tilde{\kappa}( \theta^\alpha_j) \nabla\theta^\alpha_j \cdot \nabla \psi\,\mathrm{d}x
%%- \bu^\alpha_j \theta^\alpha_j \cdot \nabla \psi   \,\mathrm{d}x
%+ \mathcal{C}[\bu^\alpha_j,\theta^\alpha_j,\psi]
%= \int_\Omega \BD\bu^\alpha_j \fp \Scr( \BD \bu^\alpha_j, \theta^\alpha_j ) \psi \,\mathrm{d}x,
%\end{equation}
\begin{equation}\label{exist:equ10}
\int_\Omega\delta\theta^\alpha_j \psi
+ \int_\Omega \tilde{\kappa}( \theta^\alpha_j) \nabla\theta^\alpha_j \cdot \nabla \psi
%- \bu^\alpha_j \theta^\alpha_j \cdot \nabla \psi   \,\mathrm{d}x
+ \mathcal{C}[\bu^\alpha_j,\theta^\alpha_j,\psi]
= \int_\Omega \BD\bu^\alpha_j \fp \Scr( \BD \bu^\alpha_j, \theta^\alpha_j ) \psi ,
\end{equation}
for every \( \bv \in \Vndiv\) and \( \psi \in U^m\).  Here \( \delta\)  denotes the backwards difference quotient of order 1, namely,
\begin{align*}
\delta\bu^\alpha_j := \frac{\bu^\alpha_j - \bu^\alpha_{j-1}}{\tau_m}, \quad \delta\theta^\alpha_j := \frac{\theta^\alpha_j - \theta^\alpha_{j-1}}{\tau}.
\end{align*}
\end{prop}

\begin{pf}
  The proof   makes use of the following corollary to Brouwer's fixed point theorem \cite[Cor.\ 1.1]{Girault1986}. The problem of finding a \( z\in X\) such that $F(z)=0$ where $F\colon X \to X$ is a function defined on a finite dimensional Hilbert space $X$ has a solution $z_*$ if there exists a  $\lambda>0$ such that $\langle F(z),z\rangle>0$ for every \( z\in X\) with $\|z\|=\lambda$. Furthermore, the solution satisfies $\|z_*\|\leq \lambda$.

  Suppose $\bu^\alpha_{j-1}\in \Vndiv$ is given, for a  $j\in\{1,\ldots,T/\tau_m\}$. For fixed $(\tilde{\bu},\tilde{\theta})\in \Vndiv\times U^m$, consider the problem of finding $\bu\in \Vndiv$ such that $F_1(\bu)=0$ where $F_1$ is defined by
%\begin{equation}\label{exist:discr_existence_eq1}
%  \langle  F_1(\bu),\bv\rangle := \int_\Omega   \frac{\bu - \bu^\alpha_{j-1}}{\tau_m} \cdot \bv
%  + \Scr ( \BD\bu , \tilde{\theta}) \fp \BD \bv
%+ \frac{1}{k}|\bu|^{r^\circledast-2}\bu \cdot \bv\,\mathrm{d}x
%+ \mathcal{B}[\tilde{\bu}, \bu, \bv]
%- \langle \bm{f}_j^m, \bv\rangle,
%%=0.
%\end{equation}
\begin{equation}\label{exist:discr_existence_eq1}
  \begin{split}
  \langle  F_1(\bu),\bv\rangle 
  & := \int_\Omega   \left[ \frac{\bu - \bu^\alpha_{j-1}}{\tau_m} \cdot \bv
  + \Scr ( \BD\bu , \tilde{\theta}) \fp \BD \bv
+ \frac{1}{k}|\bu|^{r^\circledast-2}\bu \cdot \bv\right]
\\
&
\quad+ \mathcal{B}[\tilde{\bu}, \bu, \bv]
- \langle \bm{f}_j^m, \bv\rangle,
\end{split}
\end{equation}
for $\bv\in \Vndiv$. Testing in \eqref{exist:discr_existence_eq1} with $\bv=\bu$ and using the fact that $\Scr(\BD\bu,\tilde{\theta})\fp \BD\bu \geq 0$ alongside the skew-symmetry of $\mathcal{B}$, we find that
\begin{align*}
  \langle F_1(\bu),\bu \rangle &\geq
  \frac{1}{\tau_m} \|\bu\|^2_{L^2(\Omega)}
  - \int_\Omega \frac{1}{\tau_m}\bu^\alpha_{j-1}\cdot \bu 
  - \langle \bm{f}_j^m ,\bu \rangle\\
                             &\geq
  \frac{1}{\tau_m} \|\bu\|^2_{L^2(\Omega)}
  - \varepsilon (\|\bu\|^2_{L^2(\Omega)}
  +  \|\bu\|^2_{W^{1,r}(\Omega)} ) \\
                             &\quad -  C(\varepsilon) \left(\frac{1}{\tau_m^2}\|\bu_{j-1}^\alpha\|^2_{L^2(\Omega)}
  + \|\bm{f}_j^m\|^2_{W^{-1,r'}(\Omega)}\right),
\end{align*}
where in the last step we use Young's inequality. Choosing   $\varepsilon$ sufficiently small and using the equivalence of norms in finite dimensional spaces, the aforementioned corollary to Brouwer's fixed point theorem guarantees that the solution operator $H_1\colon \Vndiv\times U^m \to \Vndiv$ given by  $H_1(\tilde{\bu},\tilde{\theta}):= \bu$ is well defined. Furthermore,  the solution satisfies $\|\bu\|_{W^{1,r}(\Omega)}\leq K_1$ where $K_1 = K_1(m,n)>0$ is independent of $\tilde{\bu}$ and $\tilde{\theta}$.

Similarly, we define $H_2\colon \Vndiv\times U^m \to U^m$  to be the solution operator associated to the function
%\begin{equation}\label{exist:discr_existence_eq2}
%\langle F_2(\theta),\psi \rangle :=
%\int_\Omega \frac{\theta - \theta^\alpha_{j-1}}{\tau_m} \psi\,\mathrm{d}x
%+ \int_\Omega \tilde{\kappa}( \theta) \nabla\theta \cdot \nabla \psi\,\mathrm{d}x
%+ \mathcal{C}[\tilde{\bu},\theta,\psi]
%- \int_\Omega \BD\tilde{\bu} : \Scr( \BD \tilde{\bu}, \tilde{\theta} ) \psi\,\mathrm{d}x.
%\end{equation}
\begin{equation}\label{exist:discr_existence_eq2}
\langle F_2(\theta),\psi \rangle :=
\int_\Omega \frac{\theta - \theta^\alpha_{j-1}}{\tau_m} \psi
+ \int_\Omega \tilde{\kappa}( \theta) \nabla\theta \cdot \nabla \psi
+ \mathcal{C}[\tilde{\bu},\theta,\psi]
- \int_\Omega \BD\tilde{\bu} : \Scr( \BD \tilde{\bu}, \tilde{\theta} ) \psi.
\end{equation}
Indeed, we set  $\theta:= H_2(\tilde{\bu},\tilde{\theta})$ if a solution \( \theta\) exists to \( F_2(\theta)  = 0 \).
Using a similar reasoning to the above, we see that, if $\|\tilde{\bu}\|_{W^{1,r}(\Omega)}\leq K_1$, the solution \( \theta\)  exists and satisfies $\|\theta\|_{H^1(\Omega)} \leq K_2$ where $K_2>0$ depends on $K_1$, $m$, and $n$, but is independent of $\tilde{\theta}$.

We observe that a fixed point of the operator $H\colon \Vndiv\times U^m \to \Vndiv\times U^m$ defined by  $H(\bu,\theta):= (H_1(\bu,\theta), H_2(\bu,\theta))$ is precisely a solution of \eqref{exist:equ9} and \eqref{exist:equ10}. By the arguments above, we see  that the operator \( H \) maps $ B^V_{K_1}\times B^U_{K_2} $ back into itself, where $B^V_{K_1}\subset \Vndiv$ and $B^U_{K_2}\subset U^m$ are the balls  of radii $K_1$ and $K_2$, respectively. Thus, if we can verify the continuity of $H_1$ and $H_2$, Brouwer's fixed point theorem   guarantees the existence of a solution, completing  the proof of the proposition.

Let us examine $H_2$ first. We take arbitrary $(\bu,\theta), (\bw,\eta)\in B^V_{K_1}\times B^U_{K_2}$ and  subtract the equations for $H_2(\bu,\theta)$ and $H_2(\bw,\eta)$. Testing in the resulting equation against the difference $\psi = H_2(\bu,\theta) - H_2(\bw,\eta)$ yields (recalling that the heat flux is monotone)
\begin{align*}
 & \frac{1}{\tau_m}\|H_2(\bu,\theta) - H_2(\bw,\eta)\|^2_{L^2(\Omega)}
  \\&\leq
  \mathcal{C}[\bw, H_2(\bw,\eta), H_2(\bu,\theta)-H_2(\bw,\eta)]
  - \mathcal{C}[\bu, H_2(\bu,\theta), H_2(\bu,\theta)-H_2(\bw,\eta)]
\\&\quad  +
  \int_\Omega\left[ (\Scr(\BD\bu,\theta)\fp \BD\bu - \Scr(\BD\bw,\eta)\fp \BD\bw)(H_2(\bu,\theta) - H_2(\bw,\eta)) \right]
  \\
%  &\quad + \frac{1}{2\tau_m} \|\theta - \eta\|_{L^2(\Omega)}^2
%  \\
  &\leq
  \mathcal{C}[\bu - \bw, H_2(\bu,\theta), H_2(\bw,\eta)]  %+ \frac{1}{2\tau_m} \|\theta - \eta\|_{L^2(\Omega)}^2
\\&\quad +
  \int_\Omega\left[  (\Scr(\BD\bu,\theta)\fp \BD\bu - \Scr(\BD\bw,\eta)\fp \BD\bw)(H_2(\bu,\theta) - H_2(\bw,\eta))\right]
  \\
  &\leq
  c_n \|\bu - \bw\|_{W^{1,r}(\Omega)}\|H_2(\bu,\theta)\|_{H^{1}(\Omega)} \|H_2(\bw,\eta)\|_{H^{1}(\Omega)}
  %+ \frac{1}{2\tau_m} \|\theta - \eta\|_{L^2(\Omega)}^2
   \\
  & \quad  + \|\Scr(\BD\bu,\theta)\fp \BD\bu - \Scr(\BD\bw,\eta)\fp \BD\bw\|_{L^2(\Omega)}
  \|H_2(\bu,\theta) - H_2(\bw,\eta)\|_{L^2(\Omega)},
\end{align*}
which, recalling the boundedness of $H_2$, implies that $H_2(\bu,\theta)\to H_2(\bw,\eta)$ as $(\bu,\theta)\to (\bw,\eta)$. We note that we rely heavily on the equivalence of norms in finite dimensional spaces in the above. In a similar fashion, reasoning for \(H_1\), we have that
\begin{align*}
 & \frac{1}{\tau_m}\|H_1(\bu,\theta) - H_1(\bw,\eta)\|^2_{L^2(\Omega)}
  \\&\leq
 %\frac{1}{2\tau_m}\|\bu - \bw\|_{L^2(\Omega)}^2 +
 \mathcal{B}[\bu - \bw, H_1(\bu,\theta), H_1(\bw,\eta)]
  \\
  & \quad +
  \int_\Omega \left[ (\Scr(\BD H_1(\bw,\eta),\eta) - \Scr(\BD H_1(\bu,\theta),\theta))\fp
  (\BD H_1(\bu,\theta) - \BD H_1(\bw,\eta))\right] 
  \\
  &\leq
  %\frac{1}{2\tau_m}\|\bu - \bw\|_{L^2(\Omega)}^2  +
  \mathcal{B}[\bu - \bw, H_1(\bu,\theta), H_1(\bw,\eta)]
  \\
  & \quad +
  \int_\Omega\left[  (\Scr(\BD H_1(\bu,\theta),\eta) - \Scr(\BD H_1(\bu,\theta),\theta))\fp
  (\BD H_1(\bu,\theta) - \BD H_1(\bw,\eta))\right] ,
%  &\leq
%  c_n \|\bu - \bw\|_{W^{1,r}(\Omega)}\|H_1(\bu,\theta)\|_{W^{1,r}(\Omega)} \|H_1(\bw,\eta)\|_{W^{1,r}(\Omega)} \\
%  & \hspace{-2.5cm} + \|\Scr(\BD H_1(\bu,\theta),\theta) - \Scr(\BD H_1(\bu,\theta),\eta)\|_{L^2(\Omega)}
%  \|H_1(\bu,\theta) - H_1(\bw,\eta)\|_{L^2(\Omega)},
\end{align*}
where in the first inequality we use that the operator $\bu \mapsto \frac{1}{k}|\bu|^{r^\circledast-2}\bu$ is monotone, since it is the derivative of a convex function. In the second inequality, we employ the monotonicity of $\Scr$. Reasoning as above, this implies that $H_1$ is continuous,  concluding the proof of the assertion.
\end{pf}

The first step towards obtaining a dissipative solution of the original problem is to take the limit in the discretisation parameter \( m \) while keeping \(\beta = (n,k) \) fixed.

\begin{lem}\label{exist:lem1}
  Let the assumptions of Proposition \ref{exist:prop1} hold and let \((\bu^{\alpha}_j, \theta^\alpha_j)_{j=1}^m \) be the sequence of solutions constructed there. \texblk{Let $\overline{\bu}^\alpha$ and $\overline{\theta}^\alpha$ denote the corresponding piecewise constant interpolants and let $\tilde{\bu}^\alpha$ the piecewise linear interpolant (recall \eqref{def:constant_interpolant} and \eqref{def:linear_interpolant}). }
%Given a sequence \((\bw_j)_{j=0}^m\), let \( \overline{\bw}\) and \( \tilde{\bw}\) denote the piecewise constant and piecewise linear interpolant, respectively.
There exist a constant \( C_1\), independent of \( \alpha\), and a constant \( C_2 = C_2(n)\), independent of \( k \) and $m$, such that
\begin{align*}
  \| &\overline{\bu}^\alpha\|_{L^\infty(0,T;L^2(\Omega))}^2
  + \tau_m \|\partial_t \tilde{\bu}^\alpha\|_{L^2(Q)}^2
  + \|\overline{\bu}^\alpha\|_{L^r(0,T;W^{1,r}(\Omega))}^r
  \\
     &\quad + \frac{1}{k}\|\overline{\bu}^\alpha\|^{r^\circledast}_{L^{r^\circledast}(Q)}
  + \|\Scr( \BD \overline{\bu}^\alpha, \overline{\theta}^\alpha) \|_{L^{r'}(Q)}^{r^\prime} \leq C_1,
\end{align*}
and
\begin{align*}
  \|\overline{\theta}^\alpha\|_{L^{\infty}(0,T;L^2(\Omega))}
  + \tau_m \|\partial_t \tilde{\theta}^\alpha\|_{L^2(Q)}^2
  + \|\nabla \theta^\alpha\|_{L^2(Q)}^2 \leq C_2.
\end{align*}
\end{lem}
\begin{pf}
Testing in (\ref{exist:equ9}) against \( \bu^\alpha_j \), we see that
\begin{align*}
&\frac{1}{2\tau_m}\Big( \|\bu^\alpha_j \|_{L^2(\Omega)}^2
- \|\bu^\alpha_{j-1}\|_{L^2(\Omega)}^2
+ \|\bu^\alpha_j  - \bu^\alpha_{j-1}\|_{L^2(\Omega)}^2 \Big)
  + \int_\Omega \Scr( \BD\bu^\alpha_j, \theta^\alpha_j)  \fp \BD \bu^\alpha_j  %\,\mathrm{d}x
\\
&\quad
+ \frac{1}{k} \|\bu^\alpha_j\|^{r^\circledast}_{L^{r^\circledast}(\Omega)}
= \langle\bm{f}^m_j , \bu^\alpha_j \rangle
\leq C(\varepsilon) \|\bm{f}^m_j \|_{W^{-1,r'}(\Omega)}^{r^\prime}
+ \varepsilon\|\bu^\alpha_j \|_{W^{1,r}(\Omega)}^r,
\end{align*}
where \( \varepsilon> 0 \) is to be fixed sufficiently small later.
By the coercivity conditions and  the Korn--Poincar\'{e} inequality, we deduce that
\begin{equation}\label{exist:equ11}
\begin{aligned}
&\|\bu^\alpha_j \|_{L^2(\Omega)}^2
 - \|\bu^\alpha_{j-1} \|_{L^2(\Omega)}^2
+ \tau_m^2 \|\delta\bu^\alpha_j \|_{L^2(\Omega)}^2
 + \tau_m \|\bu^\alpha_j \|_{W^{1,r}(\Omega)}^r
\\&\quad
+   \tau_m\|\Scr ( \BD\bu^\alpha_j , \theta^\alpha_j) \|_{L^{r'}(\Omega)}^{r^\prime}
+ \frac{\tau_m}{k} \|\bu^\alpha_j\|^{r^\circledast}_{L^{r^\circledast}(\Omega)}
\\
&\leq C(r, \varepsilon) \|\overline{\bm{f}}^m \|_{L^{r'}(t_{j-1},t_j;W^{-1,r^{\prime}}(\Omega))}^{r^\prime}
+ c \|g\|_{L^1(Q_{j-1}^j)}
+ \varepsilon\tau_m \|\bu^\alpha_j \|_{W^{1,r}(\Omega)}^r,
\end{aligned}
\end{equation}
where the function \( g\) comes from the coercivity assumption on \( \Scr\) and  we define \( Q_{j-1}^j = (t_{j-1}, t_j) \times \Omega\) for \( 1 \leq j \leq T/\tau_m\).
Choosing \( \varepsilon\) sufficiently small, we absorb the final term on the right-hand side into the left-hand side. Summing (\ref{exist:equ11}) over the indices \( j = 1,\ldots ,l \) for an \(l \in \{ 1, \ldots , T/\tau_m \}\) and maximising the resulting left-hand side, it follows that
\begin{align*}
  &\tau_m \sum_{j=1}^{T/\tau_m} \Big[\tau_m \|\delta\bu^\alpha_j \|_{L^2(\Omega)}^2
  + \|\bu^\alpha_j \|_{W^{1,r}(\Omega)}^r
+  \|\Scr(\BD\bu^\alpha_j, \theta^\alpha_j) \|_{L^{r'}(\Omega)}^{r^\prime} \Big] \\
&\quad + \max_{1\leq j \leq T/\tau_m} \|\bu^\alpha_j \|_{L^2(\Omega)}^2
+ \frac{\tau_m }{k}  \sum_{j=1}^{T/\tau_m} \|\bu^\alpha_j\|^{r^\circledast}_{L^{r^\circledast}(\Omega)}
\\&\leq C (\|\bm{f}\|^{r'}_{L^{r'}(0,T;W^{-1,r'}(\Omega))}
+ \|g\|_{L^1(Q)}
+ \|\bu_0\|^2_{L^2(\Omega)}),
\end{align*}
using the stability property
\begin{equation}\label{eq:L2_stability}
  \|P^n_V \bu_0\|_{L^2(\Omega)} \leq \|\bu_0\|_{L^2(\Omega)}.
\end{equation}
Similarly, testing in (\ref{exist:equ10}) against \( \theta^\alpha_j \) yields
\begin{align*}
&\frac{1}{2\tau_m} \Big( \|\theta^\alpha_j \|_{L^2(\Omega)}^2
- \|\theta^\alpha_{j-1}\|_{L^{2}(\Omega)}^2
+ \|\theta^\alpha_j - \theta^\alpha_{j-1}\|_{L^2(\Omega)}^2 \Big)
+ \int_\Omega \tilde{\kappa}( \theta^\alpha_j ) |\nabla \theta^\alpha_j |^2  
%- \int_\Omega \bu^\alpha_j\theta^\alpha_j  \cdot \nabla \theta^\alpha_j \,\mathrm{d}x
\\
&= \int_\Omega\Scr (\BD \bu^\alpha_j , \theta^\alpha_j ) \fp \BD \bu^\alpha_j \theta^\alpha_j
\\
& \leq C \int_\Omega\left[ ( |\BD\bu^\alpha_j |^{r-1} + 1) |\BD\bu^\alpha_j| \theta^\alpha_j \right],
%& \leq C(n)\|\Scr(\BD\bu^\alpha_j, \theta^\alpha)\|_{L^{r'}(\Omega)}
%\|\BD\bu^\alpha_j\|_{L^r(\Omega)}
%\|\theta^\alpha_j\|_{L^2(\Omega)},
\end{align*}
where we use Assumption \ref{as:const_rel} concerning the growth of \( \Scr\). However, using the fact that norms on finite-dimensional spaces are equivalent, we see that
\begin{align*}
 \int_\Omega \left[ ( |\BD\bu^\alpha_j |^{r-1} + 1) |\BD\bu^\alpha_j| \theta^\alpha_j \right]
 &\leq C \left( \|\BD\bu^\alpha_j \|^{r}_{L^{2r}(\Omega)} + \|\BD\bu^\alpha_j \|_{L^2(\Omega)} \right)  \|\theta^\alpha_j \|_{L^2(\Omega)}
 \\
% &\leq C \left( \|\BD\bu^\alpha_j \|_{L^\infty(\Omega)}^r + 1\right) \|\theta^\alpha_j \|_{L^2(\Omega)}
% \\
  &\leq C(n) \left( \|\bu^\alpha_j \|_{L^2(\Omega)}^r + \|\bu^\alpha_j\|_{L^2(\Omega)} \right) \|\theta^\alpha_j \|_{L^2(\Omega)}.
\end{align*}
However, from the above, we know for instance that
\begin{align*}
\max_{1 \leq j \leq T/\tau_m } \|\bu^\alpha_j\|_{L^2(\Omega)}^r \leq C_1^{\frac{r}{2}},
\end{align*}
where \( C_1\) is the constant from the first bound and is independent of \( \alpha\). It follows that
\begin{gather*}
 \frac{1}{2\tau_m} \Big( \|\theta^\alpha_j \|_{L^2(\Omega)}^2
- \|\theta^\alpha_{j-1}\|_{L^{2}(\Omega)}^2
+ \|\theta^\alpha_j - \theta^\alpha_{j-1}\|_{L^2(\Omega)}^2 \Big)
\\
+ \int_\Omega \tilde{\kappa}( \theta^\alpha_j ) |\nabla \theta^\alpha_j |^2 %\,\mathrm{d}x
\leq C(n) \|\theta^\alpha_j \|_{L^2(\Omega)}.
\end{gather*}
Hence, for an arbitrary \( l \in \{ 1, \dots, T/\tau_m\} \), we have
\begin{align*}
&\|\theta^\alpha_l\|_{L^2(\Omega)}^2 +  \tau_m \sum_{j=1}^l \left[ \tau_m \|\delta \theta^\alpha_j\|^2_{L^2(\Omega)}
  + \|\nabla \theta^\alpha_j\|^2_{L^2(\Omega)} \right]
  \\
&\leq C(n) \tau_m \sum_{j=1}^l \|\theta^\alpha_j\|_{L^2(\Omega)} + \|\theta^\alpha_0\|^2_{L^2(\Omega)}
  \\
  &\leq C(n) \left(  \tau_m \sum_{j=1}^l \|\theta^\alpha_j\|_{L^2(\Omega)}^2 + 1\right) .
\end{align*}
The result follows by a discrete version of Gronwall's inequality and then recalling the definition of the piecewise constant and piecewise linear interpolants.
\end{pf}

Using the {\it a priori} bounds of Lemma \ref{exist:lem1}, standard weak compactness results and Simon's lemma, the following convergence results are immediate (cf.\ \cite{Farrell2019}).

\begin{cor}\label{exist:cor1}
There exists a limiting triple \((\bu^\beta, \theta^\beta, \hat{\BS}^\beta) \) such that the following convergence results hold, up to a subsequence in \( m \) that we do not relabel:
\begin{align*}
\overline{\bu}^\alpha& \overset{\ast}{\rightharpoonup} \bu^\beta&\quad& \text{ weakly-* in }L^\infty(0, T; L^2(\Omega)^d) ;
\\
\overline{\bu}^\alpha &\rightharpoonup \bu^\beta &\quad& \text{ weakly in } L^r(0, T; W^{1,r}(\Omega)^d) ;
\\
\frac{1}{k}|\overline{\bu}^\alpha|^{r^\circledast - 2}\overline{\bu}^\alpha &\rightharpoonup \frac{1}{k}|\bu^\beta|^{r^\circledast-2}\bu^\beta &\quad& \text{ weakly in } L^{(r^\circledast)'}(Q)^d;
\\
\Scr( \BD \overline{\bu}^\alpha, \overline{\theta}^\alpha) &\rightharpoonup \hat{\BS}^\beta &\quad& \text{ weakly in }L^{r^\prime}(Q)^{d\times d};
\\
\overline{\theta}^\alpha &\overset{\ast}{\rightharpoonup} \theta^\beta
&\quad & \text{ weakly-* in } L^\infty(0, T; L^2(\Omega));
\\
\overline{\theta}^\alpha &\rightharpoonup \theta^\beta &\quad & \text{ weakly in } L^2(0, T; W^{1,2}(\Omega)) .
\end{align*}
Furthermore, we have the following strong convergence results:
\begin{align*}
\overline{\theta}^\alpha &\rightarrow\theta^\beta &\quad&\text{ strongly in }L^2(Q);
\\
\overline{\bu}^\alpha &\rightarrow\bu^\beta &\quad& \text{ strongly in }L^2(Q)^d,\\
\texblk{\overline{\bu}^\alpha} &\texblk{\rightarrow\bu^\beta} &\quad& \texblk{\text{ strongly in }L^2(0,T;W^{1,2}(\Omega)^d),}
\end{align*}
\texblk{where we use the fact that in finite dimensions all norms are equivalent (recall that at this point $n$ is fixed).}
\end{cor}
 Passing to the limit in the time discrete formulation, we see that the couple \((\bu^\beta, \theta^\beta) \) is a solution of the following semi-discrete problem. For  every \( \phi \in C^\infty_0([0, T)) \), we have the momentum balance
%\begin{equation}\label{exist:equ1}
%\begin{aligned}
%&- \int_{Q}\bu^\beta\cdot \big( \bv \partial_t\phi \big) \,\mathrm{d}x\,\mathrm{d}t - \int_\Omega P^n_V\bu_0 \cdot \big( \bv \phi(0) \big) \,\mathrm{d}x
%+ \frac{1}{k} \int_Q |\bu^\beta|^{r^\circledast-2}\bu^\beta \cdot (\bv\phi) \,\mathrm{d}x\,\mathrm{d}t
%\\&\quad
%+ \int_0^T \mathcal{B}[\bu^\beta,\bu^\beta,\bv] \phi \,\mathrm{d}t
%+ \int_Q \hat{\BS}^\beta : \BD (\bv \phi) \,\mathrm{d}x\,\mathrm{d}t
%= \int_0^T \langle \bf , \bv\rangle \phi \,\mathrm{d}t,
%\end{aligned}
%\end{equation}
\begin{equation}\label{exist:equ1}
\begin{aligned}
&- \int_{Q}\bu^\beta\cdot \big( \bv \partial_t\phi \big) 
- \int_\Omega P^n_V\bu_0 \cdot \big( \bv \phi(0) \big)
+ \frac{1}{k} \int_Q |\bu^\beta|^{r^\circledast-2}\bu^\beta \cdot (\bv\phi) 
\\&\quad
+ \int_0^T \mathcal{B}[\bu^\beta,\bu^\beta,\bv] \phi 
+ \int_Q \hat{\BS}^\beta : \BD (\bv \phi) 
= \int_0^T \langle \bf , \bv\rangle \phi ,
\end{aligned}
\end{equation}
for every \( \bv \in V^n_{\mathrm{div}}\),
and   the temperature balance
\begin{equation}\label{exist:equ2}
\begin{aligned}
&-\int_Q \theta^\beta \psi \partial_t\phi   - \int_\Omega\theta^n_0 \psi \phi   + \int_Q \tilde{\kappa} ( \theta^\beta) \nabla \theta^\beta \cdot \nabla \big( \psi \phi \big)  
\\&\quad
+ \int_0^T \mathcal{C}[\bu^\beta ,\theta^\beta, \psi] \phi  
= \int_Q \hat{\BS}^\beta : \BD\bu^\beta \psi \phi ,
\end{aligned}
\end{equation}
for every \( \psi \in W^{1,\infty}(\Omega) \).
Furthermore, as a result of the strong convergence \texblk{of the gradient of $\overline{\bu}^{\alpha}$ in $L^2(Q)^{d\times d}$}, we can identify \( \hat{\BS}^\beta  = \Scr( \BD \bu^\beta, \theta^\beta) \) a.e.\ in \( Q\).
Hence the triple \((\bu^\beta, \theta^\beta,\hat{\BS}^\beta) \) is a solution of a suitable approximation of the original problem (\ref{eq:PDE}). Now we search for appropriate uniform bounds that allow us to take the limit as \( n \rightarrow\infty \).

\begin{lem}
Let \((\bu^\beta, \theta^\beta,\hat{\BS}^\beta) \) be the solution triple constructed in Corollary \ref{exist:cor1}. There exists a positive constant \( C_1\), independent of \( \beta\), such that
\begin{equation}\label{exist:equ3}
  \|\bu^\beta\|_{L^\infty(0,T;L^2(\Omega))}^2
  + \|\bu^\beta\|_{L^r(0,T;W^{1,r}(\Omega))}^r
  +\frac{1}{k} \|\bu^\beta\|^{r^\circledast}_{L^{r^\circledast}(Q)}
  + \|\hat{\BS}^\beta\|_{L^{r'}(Q)}^{r^\prime} \leq C_1.
\end{equation}
Furthermore, there exists a positive constant $C_2$, independent of $\beta$, such that
\begin{equation}\label{exist:equ4}
  \|\bu^\beta\|_{L^{\frac{r(d+2)}{d}}(Q)} \leq C_2.
\end{equation}
\end{lem}

\begin{pf}
  The first bound (\ref{exist:equ3}) follows immediately from Lemma \ref{exist:lem1} and the weak lower semi-continuity of norms. The estimate (\ref{exist:equ4}) is a standard parabolic embedding, which is a consequence of the Gagliardo--Nirenberg inequality (see, for example, \cite[Lm.\ 7.8]{Roubicek2013}).
\end{pf}

\begin{lem}
Let \((\bu^\beta, \theta^\beta,\hat{\BS}^\beta) \) be the solution triple constructed in Corollary \ref{exist:cor1}.  There exists a constant \( c_* > 0 \), independent of \( \beta\), such that
\begin{equation}\label{exist:equ12}
\theta^\beta\geq c_*> 0 \quad\text{ a.e. in }Q.
\end{equation}
There exists a constant \( C_3>0\), independent of \( \beta\), such that
\begin{equation}\label{exist:equ5}
  \|\theta^\beta\|_{L^\infty(0,T;L^1(\Omega))}
  + \|\theta^\beta\|_{L^s(Q)}
  + \|\theta^\beta\|_{L^q(0,T;W^{1,q}(\Omega)) } \leq C_3,
\end{equation}
for   \( s\in [1,\frac{5}{3}) \) and \( q\in [1,\frac{5}{4}) \). Moreover, for sufficiently large \(p\), there exists a positive constant $C_4=C_4(k)$, independent of $n$, such that
\begin{equation}\label{exist:equ6}
  \|\partial_t \theta^\beta\|_{L^1(0,T;W^{-1,p^\prime}(\Omega)) }\leq C_4.
\end{equation}
\end{lem}

\begin{pf}
For (\ref{exist:equ12}) and \eqref{exist:equ5}, we reason as in \cite{Bulicek2009}. Although the authors there work in the setting \( r = 2\), the argument is independent of \( r\) and so can be used here.

Consequently, we see from (\ref{exist:equ2}), (\ref{exist:equ3}), (\ref{exist:equ4})  and (\ref{exist:equ5})   that (\ref{exist:equ6}) must also hold. We note that the exponent $r^\circledast$ in the penalty term was chosen to ensure that the term $\bu^\beta\cdot \nabla \theta^\beta$, which appears in the modified convective term $\mathcal{C}$, belongs to $L^{1+\delta}(Q)$, for some $\delta>0$.
\end{pf}

Similarly as in Corollary \ref{exist:cor1}, the estimates above ensure that there exists a limiting triple $(\bu^k,\theta^k,\hat{\BS}^k)$ such that the following convergences hold, up to a subsequence that we do not relabel:
\begin{align*}
\bu^\beta &\overset{\ast}{\rightharpoonup} \bu^k &\quad & \text{ weakly-* in }L^\infty(0, T; L^2(\Omega)^d) ;
\\
\bu^\beta &\rightharpoonup \bu^k &\quad & \text{ weakly in }L^r(0, T; W^{1,r}(\Omega)^d);
\\
\frac{1}{k}|\overline{\bu}^\beta|^{r^\circledast - 2}\overline{\bu}^\beta &\rightharpoonup \frac{1}{k}|\bu^k|^{r^\circledast-2}\bu^k &\quad& \text{ weakly in } L^{(r^\circledast)'}(Q)^d;
\\
\bu^\beta &\rightarrow\bu^k &\quad&\text{ strongly in } L^{q_1}(Q)^d , \, q_1\in \Big[1,r^\circledast\Big) ;
\\
\Scr^k(\BD \bu^\beta, \theta^\beta) &\rightharpoonup \hat{\BS}^k &\quad& \text{ weakly in }L^{r^\prime}(Q)^{d\times d};
\\
\theta^\beta &\rightharpoonup \theta^k &\quad & \text{ weakly in }L^{q_2}(0, T; W^{1,q_2}( \Omega)) , \, q_2\in \Big[ 1,\frac{5}{4}\Big) ;
\\
\theta^\beta&\rightarrow\theta^k &\quad & \text{ strongly in } L^{q_3}(Q),\, q_3\in \Big[ 1, \frac{5}{3}\Big) .
 \end{align*}
 {\color{black}
The strong convergence of the temperature is a consequence of the Aubin--Lions theorem and either of the estimates
\begin{equation}\label{eq:temp_time_derivative}
	\begin{aligned}
		\|		\partial_t\theta^{\beta}\|_{L^1(0,T;(W^{1,q}(\Omega))^*)}&\leq c_k 
&\qquad &\text{for }q\in (5,\infty],\\
		\|\partial_t\theta^{\beta}\|_{L^1(0,T;(W^{1,q}(\Omega))^*)}&\leq c
		&\qquad &\text{for }q\in (10,\infty],
	\end{aligned}
\end{equation}
which can be obtained from the temperature balance \eqref{exist:equ2}. Here $c_k$ is a positive constant that blows up as $k\to \infty$.
 }
 These convergence results  are sufficient to allow passage to the limit in the momentum equation. Thus we obtain
 \begin{equation}\label{eq:PDE_momentum_k}
 \int_\Omega\partial_t \bu^k \cdot \bv %\,\mathrm{d}x
 + \int_\Omega \hat{\BS}^k \fp \BD \bv%\,\mathrm{d}x
 - \int_\Omega(\bu^k \otimes\bu^k)  \fp \BD \bv %\,\mathrm{d}x
 + \frac{1}{k}\int_\Omega |\bu^k|^{r^\circledast-2}\bu^k\cdot \bv %\,\mathrm{d}x
 = \langle \bf, \bv\rangle,
 \end{equation}
 for every \( \bv\in C^\infty_{0,\diver}(\Omega)^d\) and a.e.\ $t\in (0,T)$. The convective term can be now written in its original form since $\diver \bu^k =0$ pointwise. Furthermore, we claim that   $\hat{\BS}^k = \Scr(\BD\bu^k,\theta^k)$ a.e.\ in $Q$. Indeed, since $k$ is fixed, the velocity $\bu^k$ is an admissible test function in \eqref{eq:PDE_momentum_k}. Therefore we have an energy identity available (cf.\ \cite[Eq.\ 4.103]{Sueli2018}). This makes it straightforward to prove that
 \begin{equation}\label{eq:limsup_n}
   \limsup_{n\to\infty} \int_Q \Scr(\BD\bu^\beta,\theta^\beta)\fp \BD\bu^\beta 
   \leq \int_Q \hat{\BS}^k\fp \BD\bu^k .
 \end{equation}
Recalling the growth condition \eqref{eq:growth}, we observe that the dominated convergence theorem  implies that, for an arbitrary $\Btau\in L^{r}(Q)^{d\times d}$, we have
 \begin{equation}\label{eq:minty_n_1}
   \Scr(\Btau, \theta^\beta) \to \Scr(\Btau, \theta^k) \qquad
   \text{strongly in }L^{r'}(\Omega)^{d\times d},
 \end{equation}
 as $n\to \infty$. Combining the monotonicity property \eqref{eq:monotonicity} of $\Scr$ with \eqref{eq:limsup_n} and \eqref{eq:minty_n_1} yields, for an arbitrary $\Btau\in L^r(Q)^{d\times d}$,
\begin{align*}
  0 &\leq \limsup_{n\to\infty} \int_Q \left[  (\Scr(\BD\bu^\beta,\theta^\beta) - \Scr(\Btau,\theta^\beta ))\fp (\BD\bu^\beta - \Btau) \right]
  \\
    &\leq \int_Q \left[ (\hat{\BS}^k - \Scr(\Btau,\theta^k)) \fp (\BD\bu^k - \Btau)\right].
\end{align*}
Choosing $\Btau = \BD\bu^k \pm \varepsilon\Bsigma$ for  an arbitrary $\Bsigma\in C^\infty_0(Q)^{d\times d}$ and letting $\varepsilon \to 0$ concludes the proof of the claim.

In order to pass to the limit in the temperature equation, we need to investigate the convergence of $\Scr(\BD\bu^\beta,\theta^\beta)\fp \BD\bu^\beta$ in $L^1(Q)$. Firstly, from the monotonicity of $\Scr$, we   see that
\begin{equation*}
  \int_Q \hat{\BS}^k\fp \BD\bu^k 
  \leq \liminf_{n\to \infty} \int_Q \Scr(\BD\bu^\beta,\theta^\beta )\fp \BD\bu^\beta ,
\end{equation*}
and so, by \eqref{eq:limsup_n}, the equality actually holds. In turn, this implies that
\begin{equation*}
  (\Scr(\BD\bu^\beta,\theta^\beta) - \Scr(\BD\bu^k,\theta^k)) \fp
  (\BD\bu^\beta - \BD\bu^k) \to 0
  \quad \text{strongly in } L^1(Q),
\end{equation*}
noting that the sequence on the left-hand side is non-negative.
Writing
\begin{align*}
  \Scr(\BD\bu^\beta,\theta^\beta) \fp \BD\bu^\beta &=
  \Scr(\BD\bu^\beta,\theta^\beta)\fp \BD\bu^k
  + \Scr(\BD\bu^k,\theta ^k)\fp (\BD\bu^\beta - \BD\bu^k) \\
                                                   &\quad +
  (\Scr(\BD\bu^\beta,\theta^\beta) - \Scr(\BD\bu^k,\theta^k)) \fp
  (\BD\bu^\beta - \BD\bu^k),
\end{align*}
immediately yields that $\Scr(\BD\bu^\beta,\theta^\beta)\fp \BD\bu^\beta \rightharpoonup \Scr(\BD\bu^k,\theta^k) : \BD \bu^k$ weakly in $L^1(Q)$ as $n\to\infty$. Hence, the limiting functions satisfy the temperature balance
\begin{equation}\label{eq:PDE_temperature_k}
\int_\Omega \partial_t \theta^k \psi %\,\mathrm{d}x
+ \int_\Omega \tilde{\kappa}(\theta^k)\nabla \theta^k \cdot \nabla \psi %\,\mathrm{d}x
- \int_\Omega \bu^k\theta^k \cdot \nabla \psi %\,\mathrm{d}x
= \int_\Omega \Scr(\BD\bu^k,\theta^k)\fp \BD\bu^k\psi %\,\mathrm{d}x,
\end{equation}
for every $\psi \in W^{1,\infty}(\Omega)$ and a.e. \( t \in (0, T) \).  \texblk{Here we note that the initial conditions are satisfied in the following sense (see e.g.\ \cite[Lm.\ 8]{Maringova2018}):}
{\color{black}
\begin{equation}\label{eq:initial_cond_k}
\lim_{t\to 0} \|\bu^k(t,\cdot) - \bu_0\|_{L^2(\Omega)} = 0
\qquad 
\esslim_{t\to 0} \|\theta^k(t,\cdot) - \theta_0\|_{L^1(\Omega)} = 0.
\end{equation}
}
Using the weak lower semi-continuity of norms, we obtain the following estimates:
\begin{equation}\label{eq:estimates_k}
\begin{split}
  \|&\bu^k\|_{L^\infty(0,T;L^2(\Omega))}
  + \|\bu^k\|_{L^r(0,T;W^{1,r}(\Omega))}
  +\frac{1}{k} \|\bu^k\|^{r^\circledast}_{L^{r^\circledast}(Q)}
  + \|\Scr(\BD\bu^k,\theta^k)\|_{L^{r'}(Q)}  \\
    & + \|\bu^k\|_{L^{\frac{r(d+2)}{d}}(Q)}
  + \|\theta^k\|_{L^\infty(0,T;L^1(\Omega))}
  + \|\theta^k\|_{L^s(Q)}
  + \|\theta^k\|_{L^q(0,T;W^{1,q}(\Omega))}
  \leq C_1,
\end{split}
\end{equation}
for arbitrary \( s\in [1,\frac{5}{3}) \) and \( q\in [1,\frac{5}{4}) \), and a constant $C_1$ that is independent of $k$. \texblk{Note also that the almost everywhere convergence of $\theta^\beta$ with \eqref{exist:equ12} implies that $\theta^k\geq c_*>0$.} It follows that, up to a subsequence in $k$ that we do not relabel, the following convergence results hold:
\begin{align*}
\bu^k &\overset{\ast}{\rightharpoonup}\bu  &\quad& \text{ weakly-* in } L^\infty(0, T; L^2(\Omega)^d) ;
\\
\bu^k &\rightharpoonup \bu &\quad& \text{ weakly in }L^r(0, T; W^{1,r}(\Omega)^d) ;
\\
\bu^k &\rightarrow\bu &\quad& \text{ strongly in } L^{q_1}(Q)^d, \, q_1\in \Big[1,\frac{r(d+2)}{d}\Big) ;
\\
\bu^k(s,\cdot) &\rightarrow\bu(s,\cdot) &\quad& \text{ strongly in } L^{2}(\Omega)^d, \, \text{for a.e. }s\in (0,T);
\\
\frac{1}{k}|\bu^k|^{r^\circledast -2}\bu^k &\rightarrow\bm{0} &\quad& \text{ strongly in } L^{1}(Q)^d;
\\
%\partial_t\bu^k &\rightharpoonup \partial_t \bu&\quad & \text{ weakly in } L^{\check{r}}(0, T; W^{-1, \check{r}}( \Omega)^d) ;
%\\
\Scr^k(\BD \bu^k, \theta^k) &\rightharpoonup \BS &\quad & \text{ weakly in }L^{r^\prime}(Q)^{d\times d};
\\
\theta^k & \rightharpoonup \theta &\quad & \text{ weakly in } L^s(0, T; W^{1,q_2}(\Omega)) , \, q_2\in \Big[ 1, \frac{5}{4}\Big) ;
\\
\theta^k & \rightarrow\theta &\quad & \text{ strongly in } L^{q_3}(Q), \, q_3\in \Big[ 1, \frac{5}{3}\Big);
\\
\theta^k(s,\cdot) & \rightarrow\theta(s,\cdot) &\quad & \text{ strongly in } L^{1}(\Omega), \,\text{for a.e. }s\in (0,T).
\end{align*}
The above convergences do not suffice to pass to the limit in the temperature equation. Hence, at this stage we turn to the entropy balance. The following lemma states the a priori estimates that are satisfied by the entropy, which we use when passing to the limit in \( k \).

\begin{lem}\label{exist:lem2}
  Let the assumptions of Proposition \ref{exist:prop1} hold and let \(( \bu^k, \theta^k, \hat{\BS}^k) \) be the   limiting solution of \eqref{eq:PDE_momentum_k}, \eqref{eq:PDE_temperature_k} that is constructed by taking the limit as \( m \rightarrow\infty \), then \( n \rightarrow\infty \) for the solution triple from Proposition \ref{exist:prop1}. Define the entropy by \( S^k = \log\theta^k\) (without loss of generality setting $c_v \equiv 1$). There exists a constant \( C_5 \),  independent of \( k\),  such that
\begin{equation}\label{eq:estimates_k_S}
  \|S^k\|_{L^2(0,T;W^{1,2}(\Omega)) }
  + \|S^k\|_{L^\infty(0,T;L^q(\Omega))}\leq C_5,
\end{equation}
for arbitrary \( q\in [1,\infty) \).
\end{lem}
\begin{pf}
 We notice that,  for any \( \lambda\in (-1,0)\), the function \(f(x) = \frac{\log^2 x}{x^{\lambda+1}}\) is bounded on \( [1,\infty) \). Assuming that $\theta^k \geq 1$, it follows that \( \log^2\theta^k\leq c \theta^{k^{\lambda+1}} \) for a constant \( c>0 \).
  On the other hand, now suppose that \( c_*<1 \) and consider when \( \theta^k \in [c_*, 1) \).
 We notice that   \( \log^2\theta^k \leq \log^2 c_* \) and so we can bound
\begin{align*}
\frac{\log^2\theta^k}{\theta^{k^{\lambda + 1}}} \leq \frac{\log^2 c_*}{c_*^{\lambda+1}}
	 < \infty.
\end{align*}
Consequently $\log^2\theta^k \leq c {\theta^k}^{\lambda+1}$ a.e.\ in $Q$ for some positive constant $c$ that is independent of \( k \). Furthermore, from the inequality $0< c_* \leq \theta^k$, we   have that \( \theta^{k^{-2}} \leq c_*^{- 1-\lambda} \theta^{k^{\lambda -1}} \).   Hence, we deduce that
\begin{align*}
  \|S^k\|_{L^2(0,T;W^{1,2}(\Omega)) }
  = \|\log\theta^k\|_{L^2(0,T;W^{1,2}(\Omega)) }
  \leq c\|\theta^{k^{\frac{1 + \lambda}{2}}}\|_{L^2(0,T;W^{1,2}(\Omega)) } \leq c,
\end{align*}
\texblk{where we used the fact that ${\theta^k}^{\frac{1+\lambda}{2}}$ is bounded uniformly in $L^2(0,T;W^{1,2}(\Omega))$, for $\lambda\in (-1,0)$. For a proof we refer the reader to \cite{Bulicek2009}. Employing the same argument now with the function $f(x)= \frac{\log^q x}{x}$},  we deduce the latter bound in the statement of the lemma as required.
\end{pf}

As a consequence of Lemma \ref{exist:lem2}, up to a subsequence in $k$ that we do not relabel, we have the convergence results
\begin{align*}
S^k&\rightharpoonup S &\quad & \text{ weakly in } L^2(0, T; W^{1,2}(\Omega)) ,
\\
S^k &\overset{\ast}{\rightharpoonup} S &\quad & \text{ weakly-* in }L^\infty(0, T; L^{q_5}(\Omega)), \, {q_5}\in [1,\infty).
 \end{align*}
Additionally, the almost everywhere convergence of $\theta^k$ in \( Q\) allows us to identify $S = \log\theta$. We   now focus on proving that the limiting functions constitute a dissipative weak solution of \eqref{eq:PDE} in the sense of Definition \ref{def:dissipative_weak_sols}. 

\begin{thm}\label{thm:existence_end}
  Let the assumptions of Proposition \ref{exist:prop1} hold. Then, the quartet $(\bu,\theta,\BS,S)$ constructed above is a dissipative weak solution of \eqref{eq:PDE} in the sense of Definition \ref{def:dissipative_weak_sols}.
\end{thm}

\begin{pf}
  The estimates in \eqref{eq:estimates_k} and the convergences they induce make it straightforward to pass to the limit in the momentum equation and deduce that the couple $(\bu,\BS)$ satisfies \eqref{eq:weakPDE_momentum}.

  The constitutive relation \eqref{eq:weakPDE_CR} can be identified in the same manner as before, assuming that we can prove an  estimate that is analogous to \eqref{eq:limsup_n}. Since   the velocity $\bu$ is not an admissible test function in the balance of momentum, there is no energy identity available. Hence obtaining such an estimate is not straightforward. Thankfully, this difficulty can be dealt with by testing   with a Lipschitz truncation of the error $\bu^k - \bu$. In doing so, it is possible to prove the existence of a nonincreasing sequence of sets $E_j\subset Q$ such that $|E_j|\to 0$ as $j\to \infty$ and
  \begin{equation}\label{eq:limsup_k}
    \limsup_{k\to\infty} \int_{Q\setminus E_j} \Scr(\BD\bu^k,\theta^k) \fp \BD\bu^k 
    = \int_{Q\setminus E_j} \BS\fp \BD\bu 
    \qquad \text{for any }j\in \mathbb{N}.
  \end{equation}
\texblk{The general idea is that in the regions where the error is Lipschitz, an energy identity is available (since one can test the momentum equation) and an inequality analogous to \eqref{eq:limsup_n} can be obtained. It is crucial as well that the size of the ``bad set'' can be controlled.}  The details of this argument can be found, for example,  in \cite[Thm.\ 3.3]{Blechta2019}. This   implies that $\BS = \Scr(\BD\bu,\theta)$ a.e. in $Q\setminus E_j$ and that $\hat{\BS}^k\fp \BD\bu^k \rightharpoonup \BS\fp \BD\bu$ weakly in $L^1(Q\setminus E_j)$ for any $j\in \mathbb{N}$. The measure of the sets $E_j$ vanishes as \( j \rightarrow\infty\) and so   identification of the constitutive relation \eqref{eq:weakPDE_CR} follows.

Next, we turn to the entropy balance. Testing the temperature equation \eqref{eq:PDE_temperature_k} with $\psi/\theta^k$, where $\psi\in C_0^\infty([0,T);C^1(\overline{\Omega}))$ is an arbitrary function such  that $\psi\geq 0$, we obtain the following equation for the entropy $S^k$:
\begin{equation}\label{equ:entropy_kapprox}
\begin{aligned}
	&-\int_Q S^k \partial_t \psi  %\,\mathrm{d}t
	- \int_\Omega \psi(0) \log\theta_0 %\,\mathrm{d}x
	-\int_Q S^k\bu^k\cdot \nabla \psi %\,\mathrm{d}x\,\mathrm{d}t
+ \int_Q \frac{\tilde{\kappa}(\theta^k)\nabla\theta^k}{{\theta^k}}\cdot \nabla \psi %\,\mathrm{d}x\,\mathrm{d}t
\\
&\quad = \int_Q \frac{\tilde{\kappa}(\theta^k)|\nabla\theta^k|^2}{{\theta^k}^2}\psi %\,\mathrm{d}x\,\mathrm{d}t
+ \int_Q \frac{\Scr(\BD\bu^k,\theta^k)\fp \BD\bu^k}{\theta^k}\psi %\,\mathrm{d}x\,\mathrm{d}t
\geq 0
,
\end{aligned}
\end{equation}
for every \( \psi \in C_0^\infty([0,T);C^1(\overline{\Omega}))\) such that  \(\psi\geq 0\).
Now, cf.\ \cite{Maringova2018}, we claim that
\begin{equation}
  \liminf_{k\to \infty} \int_Q \texblk{\frac{\hat{\BS}^k\fp \BD\bu^k}{\theta^k} \psi}  \geq \int_Q \texblk{\frac{\BS\fp \BD\bu}{\theta}\psi} ,
\end{equation}
for any non-negative function \texblk{$\psi \in C_0^\infty([0,T);C^1(\overline{\Omega}))$}. Since \texblk{$\frac{\Scr(\BD\bu,\theta)\fp \BD\bu}{\theta}\psi$} is integrable, for every $\varepsilon>0$ there exists a $\delta>0$ such that, for any $E\subset Q$ with $|E|\leq \delta$, we have
\begin{equation*}
\int_{\texblk{E}} \texblk{\frac{\Scr(\BD\bu,\theta)\fp \BD\bu}{\theta}} \psi  \leq \varepsilon.
\end{equation*}
Noting that \texblk{the integrand} is non-negative and choosing $j\in\mathbb{N}$ sufficiently  large  so that $|E_j|\leq \delta $ for the sets $E_j$ described in \eqref{eq:limsup_k}, we have
\begin{gather*}
  \liminf_{k\to\infty} \int_Q \texblk{\frac{\hat{\BS}^k\fp \BD\bu^k}{\theta^k}} \psi %\,\mathrm{d}x\,\mathrm{d}t
  \geq \liminf_{k\to\infty} \int_{Q\setminus E_j} \texblk{\frac{\hat{\BS}^k\fp \BD\bu^k}{\theta^k}} \psi %\,\mathrm{d}x\,\mathrm{d}t
  \\
  \texblk{ =  \liminf_{k\to\infty} \left[ \int_{Q\setminus E_j} \frac{\hat{\BS}^k\fp \BD\bu^k}{\theta}\psi
      + \int_{Q\setminus E_j} \hat{\BS}^k\fp\BD\bu^k\left(\frac{1}{\theta^k}-\frac{1}{\theta}\right)\psi
    \right]
  }
  \\
\texblk{\geq} \int_{Q\setminus E_j} \texblk{\frac{\BS\fp \BD\bu}{\theta}} \psi %\,\mathrm{d}x\,\mathrm{d}t
\geq \int_{Q} \texblk{\frac{\BS\fp \BD\bu}{\theta}} \psi%\,\mathrm{d}x\,\mathrm{d}t  
- \varepsilon,
\end{gather*}
\texblk{where we used the fact that $\hat{\BS}^k\fp \BD\bu^k\wconv \BS\fp\BD\bu$ weakly in $L^1(Q\setminus E_j)$ for any $j\in\mathbb{N}$, and that $\theta^k$ converges pointwise almost everywhere to $\theta$.} Combining this with the weak lower semi-continuity of the $L^2$-norm alongside the \texblk{weak convergence in $L^2(Q)$ of the sequence $\{\sqrt{\tilde{\kappa}(\theta^k)}\nabla S^k\}_k$}, we can take the limit in the equation (\ref{equ:entropy_kapprox}) for $S^k$   and obtain the entropy inequality \eqref{eq:weakPDE_entropy}. The assumption $r>\frac{2d}{d+2}$ guarantees the compactness of $S^k\bu^k$ in the advective term.

Finally, we consider  the balance of total energy. First, we note that, since the velocity $\bu^k$ is an admissible test function in \eqref{eq:PDE_momentum_k}, the following energy identity holds for a.e. $\tau \in (0,T)$ \cite[Eq.\ 4.103]{Sueli2018}:
\begin{equation}\label{eq:energy_identity_k}
  \frac{1}{2}\|\bu^k(\tau,\cdot)\|^2_{L^2(\Omega)}
+ \int_{Q_\tau} \BS^k\fp \BD\bu^k %\,\mathrm{d}x\,\mathrm{d}t
+ \frac{1}{k}\|\bu^k\|^{r^\circledast}_{L^{r^\circledast}(Q_\tau)}
= \int_0^\tau \langle \bm{f}, \bu^k\rangle %\,\mathrm{d}t
+  \frac{1}{2}\|\bu_0\|^2_{L^2(\Omega)}.
\end{equation}
Testing the balance of temperature \eqref{eq:PDE_temperature_k} with the approximate indicator function $\psi^j$ of the interval $(0,\tau)$ and letting $j\to \infty$, we add the result to \eqref{eq:energy_identity_k} to obtain
\begin{equation*}
\Big[ \int_\Omega \left(  \frac{|\bu^k(t,\cdot) |^2}{2}
+ \theta^k(t,\cdot)   \right) \Big]_{t = 0}^{t = \tau}
+ \frac{1}{k}\|\bu^k\|^{r^\circledast}_{L^{r^\circledast}(Q_\tau)}
= \int_0^\tau \langle\bf, \bu^k\rangle \,\mathrm{d}t,
\end{equation*}
which, taking \( k \rightarrow\infty\),   implies the balance of total energy \eqref{eq:weakPDE_energy}. \texblk{At this level, the analogous statement to \eqref{eq:initial_cond_k} also holds true and the proof is identical to the one from \cite{Maringova2018}. As a consequence, from Vitali's convergence theorem we also have that} 
{\color{black}
\begin{equation*}
\esslim_{t\to 0}\|S(t,\cdot) - \log\theta_0(\cdot)\|_{L^q(\Omega)} =0 
\qquad q\in [1,\infty).
\end{equation*}
}
\end{pf}
%\AG{Something I hadn't thought about...   in what sense does the entropy satisfy the initial condition? We know that $\esslim_{t\to 0}\|\theta(t)-\theta_0\|_{L^1}=0$; what does this tell us about $\log S$?}
\begin{rmk}\label{rem:implicitCR}
  So far we have assumed that the constitutive relation is of the specific form $\BS = \Scr(\BD\bu,\theta)$, but the  approach presented here can be almost identically applied to models of the form $\BD\bu = \Dcr(\BS,\theta)$, which include for instance Glen's model for ice dynamics \cite{Glen1955}. In fact, this is true also for implicit models in which the constitutive relation is written as
\begin{equation}\label{eq:implicit_CR}
  \BG(\cdot,\BS,\BD\bu) = \bm{0}\qquad
  \text{a.e. in }Q,
\end{equation}
where $\BG\colon Q\times \Rds\times \Rds \to \Rds$ is a function defining a \emph{maximal monotone $r$-graph} (see \cite{Blechta2019} for an in-depth discussion of these models). An important example of a model that can be described in such manner  is the Herschel--Bulkley model for viscoplastic flow:
\begin{equation}\label{eq:HerschelBulkley}
  \renewcommand{\arraystretch}{1.2}
  \left\{
    \begin{array}{ccc}
      |\BS|\leq \tau_* & \Longleftrightarrow & \BD\bu = \bm{0},\\
      |\BS| > \tau_* & \Longleftrightarrow & \BS = K|\BD\bu|^{r-2}\BD\bu+ \displaystyle\frac{\tau_*}{|\BD\bu|}\BD\bu,
    \end{array}
    \right.
\end{equation}
where $\tau_*$, $K>0$ are parameters.  We note that this model can also be  described using the implicit function
\begin{equation}\label{eq:HB_implicit}
  \BG(\BS,\BD\bu) := (|\BS|-\tau_*)^+\BS
  - K|\BD\bu|^{r-2}(\tau_* + (|\BS|- \tau_*)^+)\BD\bu.
\end{equation}
For the model with $r=2$, the model with temperature dependent activation parameters can also be included (cf.\ \cite{Maringova2018}). When tackling the question of existence of dissipative weak solutions for implicit models one could employ a finite element formulation including the stress $\BS$ as an unknown, in an analogous way to the formulations analysed in \cite{Farrell2020b,Farrell2019}.
\end{rmk}

\begin{rmk}
  As mentioned in the introduction, the approximation scheme  introduced in Proposition \ref{exist:prop1} is better suited to numerical analysis than the one proposed in \cite{Bulicek2009}, since we do not require the computation of a Helmholtz decomposition nor a quasi-compressible approximation. The admissibility problem of the velocity is instead dealt with by means of a penalty term.
\end{rmk}

\begin{rmk}
  The arguments presented in this section can be applied almost verbatim to the problem with Navier's slip boundary conditions, in which case one obtains also the existence of an integrable pressure $p\in L^1(Q)$. Assuming that $r>\frac{3d}{d+2}$, the results from \cite{Bulicek2009a} guarantee even that a weak version of the energy balance \eqref{eq:PDE_total_energy} holds. The results from this paper then constitute an extension to the regime $r\in (\frac{2d}{d+2}, \frac{3d}{d+2}]$ for problems with such boundary conditions.
\end{rmk}

\begin{rmk}\label{rem:3lev_approx}
  The proof of Theorem \ref{thm:existence_end} is based on a 3-level approximation scheme that in practice could be very likely simplified. The reason for considering two discretisation indices $m$ and $n$ is to simplify the argument for obtaining the positivity of the temperature \eqref{exist:equ12} and the estimates \eqref{exist:equ5}, since otherwise these properties would have to be obtained at the finite element level. The index $k$ associated to the penalty term can very likely also be avoided, but a discrete version of a parabolic Lipschitz truncation would be needed, which although very plausible, is not available at the time of this publication (a steady version was developed in \cite{Diening:2013}).
\end{rmk}

\section{Weak-strong uniqueness}\label{sec:WSU}
In this section, we prove a weak-strong uniqueness result for weak solutions of (\ref{eq:PDE}). To that end, we introduce the following assumptions on the heat conductivity function \( \tilde{\kappa}\) and the function \( \Scr\) that defines the constitutive relation between \( \BS\) and \( \BD\bu \).

\begin{ass}\label{as:WSU_assumptions}
The functions $\tilde{\kappa}\colon \R \to \R$ and $\Scr\colon \Rds\times \R \to \Rds$ defining the constitutive relations are continuous on their domain and \( \Scr\)  satisfies the coercivity, growth and compatibility conditions from Assumption \ref{as:const_rel}. 
The following further properties are also satisfied. 
\begin{itemize}[leftmargin = 0.8cm]
  \item The function $\tilde{\kappa}$ is locally Lipschitz continuous and there exist constants $c_1$, $c_2$ such that
\begin{align*}
0 < c_1 \leq \tilde{\kappa}(s) \leq c_2(s^{\frac{1}{2}} + 1).
\end{align*}
\item For every fixed $s\in \R$ and for every $\Btau_1,\Btau_2\in \Rds$, the function $\Scr$ satisfies the strong monotonicity condition
  \begin{equation}\label{eq:strong_monotonicity}
    (\Scr(\Btau_1,s) - \Scr(\Btau_2,s))\fp (\Btau_1 - \Btau_2)
    \geq c(|\Scr(\Btau_1,s) - \Scr(\Btau_2,s)|^2 + |\Btau_1 - \Btau_2|^2),
  \end{equation}
  for some positive constant $c>0$.
\item For every $\delta>0$ and $R>0$, there exists a constant $C=C(\delta,R)$ such that
  \begin{equation}\label{eq:loc_Lipschitz_S}
|\Scr(\Btau, \eta_1) - \Scr(\Btau, \eta_2)| \leq C|\eta_1 - \eta_2|,
  \end{equation}
  for every $\eta_1,\eta_2\in [\delta,\delta^{-1}]$ and $\Btau\in \Rds$ with $|\Btau|\leq R$.
\end{itemize}
\end{ass}

This assumption can be suitably modified for the problem with constitutive relations of the kind $\BD\bu = \Dcr(\BS,\theta)$ or implicit relations. We further assume for simplicity that the body force is not present, that is, \( \bf = \mathbf{0}\). We focus on proving the following result.

\begin{thm}\label{thm:wsu}
Suppose that Assumption \ref{as:WSU_assumptions} holds. Assume  that the data \( \bu_0 \in L^2_{\mathrm{div}}(\Omega)^d\) and \( \theta_0\in L^1(\Omega) \) are given such that \( \theta_0 \geq c_* >0 \) for a constant \( c_*\).  Let \((\bu, \theta) \) be a weak solution of (\ref{eq:PDE}) as constructed in Theorem \ref{thm:existence_end}. Let \((\tilde{\bu}, \tilde{\theta})  \) be a classical solution of (\ref{eq:PDE}) with initial data \((\tilde{\bu}_0, \tilde{\theta}_0) \). The following comparison inequality holds for a constant \texblk{$C=C(\|\tilde{\theta}\|_{W^{1,\infty}(Q)}, \|\tilde{\bu}\|_{W^{1,\infty}(\Omega)})$} depending only on  the smooth solution, the lower bound \( c_* \) of \( \theta\) and the initial data \((\bu_0, \theta_0) \), for a.e. time \( t\in (0, T) \):
\begin{align*}
\int_\Omega \mathcal{E}(\bu(t), \theta(t) \,|\, \tilde{\bu}(t), \tilde{\theta}(t))   \leq \Big( \int_\Omega\mathcal{E}(\bu_0, \theta_0\,|\, \tilde{\bu}_0, \tilde{\theta}_0)   \Big)\mathrm{e}^{Ct},
\end{align*}
where we define the relative energy \( \mathcal{E}\) by
\begin{align*}
\mathcal{E}(\bu, \theta\,|\, \tilde{\bu}, \tilde{\theta}) = \frac{1}{2}|\bu - \tilde{\bu}|^2 + (\theta - \tilde{\theta}) + \tilde{\theta}( \log\tilde{\theta} - \log
{\theta}) .
\end{align*}
In particular, if the classical solution emanates from the same data as the weak solution,   weak-strong uniqueness holds so
\begin{align*}
\bu = \tilde{\bu}, \quad \theta = \tilde{\theta} \quad \text{ in }Q.
\end{align*}
\end{thm}

A major issue in the analysis is that \( \theta\) is possibly unbounded above. To overcome this,  we use the idea of essential and residual parts of a function, depending on the value of \( \theta\). Similar ideas are used in \cite{Feireisl2012}.
We fix \( \delta>0 \) sufficiently small such that \( \tilde{\theta} \in [2\delta, (2\delta)^{-1}]\) in \( Q\). Let \( \psi = \psi_\delta\in C^\infty_c([0, \infty)) \) be such that \( \psi \in [0, 1]\) with \( \psi = 1\) on \( [2\delta, (2\delta)^{-1}]\) and \( \psi = 0 \) outside of \( [\delta, \delta^{-1}]\). Given a \texblk{measurable} function \texblk{$h$}, we define the essential part \( h_{ess}\) and residual  part \( h_{res}\) of \( h \) by
{\color{black}
\begin{align*}
h_{ess}:= h\psi(\theta), \quad h_{res} := h - h_{ess} = (1 - \psi(\theta)) h.
\end{align*}
}
%\AG{I rewrote this slightly because I don't think it was correct...}
%\VP{I think it was fine as it was but maybe misleading in that it suggested that \( h \) depended only the value of \( 
%\theta\) which is generally not true in our case. Either way, let's keep it as it is written here.}
With this in mind, we recall the following useful fact:  there exists a constant \( C\), depending only on the smooth solution \((\tilde{\bu}, \tilde{\theta})\) (and thus \( \delta\)), such that
\begin{equation}\label{wsu:equ11}
\mathcal{E}(\bu, \theta \,|\, \tilde{\bu}, \tilde{\theta}) \geq C\Big\{  \big[  |\bu - \tilde{\bu} |^2 + |\theta - \tilde{\theta}|^2  \big]_{ess} + \big[ 1  +  |\log\theta|  + \theta \big]_{res}\Big\},
\end{equation}
pointwise a.e. in \( Q\) \cite{Brezina2018}. 

%\AG{Do we ever use the $[|\bu|^2]_{res}$? I think not... we could just drop it}
%\VP{No we don't use it as far as I'm aware so we can drop it if you want.}

With this, we are now ready to prove Theorem \ref{thm:wsu}. We start by using the properties of \((\bu, \theta) \) being a weak solution of (\ref{eq:PDE}). In the following $\BS$ and $\tilde{\BS}$ will denote the stresses corresponding to the weak solution and strong solution, respectively, so  $\BS = \Scr(\BD\bu,\theta)$ and $\tilde{\BS}=\Scr(\BD\tilde{\bu},\tilde{\theta})$.
Testing against \( \bv = \tilde{\bu}\) in (\ref{eq:weakPDE_momentum}), we use the pointwise divergence free property of \( \tilde{\bu}\) and weak divergence free property of \( \bu\) to see that
\begin{align*}
&\int_{Q_\tau} \bu \cdot \partial_t \tilde{\bu} %\,\mathrm{d}x\,\mathrm{d}t 
- \Big[ \int_\Omega \bu \cdot \tilde{\bu}  \Big]_{t = 0}^{t= \tau}
= \int_{Q_\tau} \BS \fp \BD \tilde{\bu} 
- (\bu - \tilde{\bu} ) \otimes (\bu - \tilde{\bu}) \fp \BD \tilde{\bu},
\end{align*}
where we denote \( Q_{\tau} = (0, \tau ) \times \Omega\).
It follows that
\begin{align*}
\Big[ \int_\Omega \frac{|\bu - \tilde{\bu}|^2 }{2}  \Big]_{t = 0}^{t = \tau}
&= \Big[ \int_\Omega \frac{|\bu |^2}{2}  \Big]_{t =0}^{t = \tau} + \int_{Q_\tau}   \tilde{\bu} \cdot \partial_t \tilde{\bu}   - \Big[ \int_\Omega\bu \cdot \tilde{\bu}   \Big]_{t = 0}^{t = \tau}
\\
&= \Big[ \int_\Omega \frac{|\bu |^2}{2}   \Big]_{t =0}^{t = \tau} - \int_{Q_\tau} ( \bu -   \tilde{\bu})  \cdot \partial_t \tilde{\bu}   \\&\quad
+ \int_{Q_\tau} \left[ \BS : \BD \tilde{\bu} - (\bu - \tilde{\bu}) \otimes (\bu - \tilde{\bu}) : \BD \tilde{\bu} \right].
\end{align*}
Using the total energy balance for the weak solution, we replace the first term on the right-hand side to deduce that
\begin{equation}\label{wsu:equ1}
\begin{aligned}
&\Big[ \int_\Omega \left( \frac{|\bu - \tilde{\bu}|^2}{2} + \theta\right) \Big]_{t = 0}^{t = \tau} - \int_{Q_\tau} \BS : \BD \tilde{\bu} 
\\
&\leq  -\int_{Q_\tau}  \left[ (\bu - \tilde{\bu}) \cdot \partial_t \tilde{\bu} + (\bu - \tilde{\bu}) \otimes (\bu - \tilde{\bu}) :\BD \tilde{\bu} \right] .
\end{aligned}
\end{equation}
Using the fact that \((\tilde{\bu}, \tilde{\theta}) \) is a classical solution and, in particular, satisfies (\ref{eq:PDE_momentum}) pointwise, we rewrite the first term on the right-hand side of (\ref{wsu:equ1}) as
\begin{equation}\label{wsu:equ2}
\int_{Q_\tau}  (\bu - \tilde{\bu}) \cdot \partial_t \tilde{\bu} 
= \int_{Q_\tau} ( \bu - \tilde{\bu}) \cdot \diver (\tilde{\BS} - \tilde{\bu}\otimes \tilde{\bu}) .
\end{equation}
Applying the appropriate divergence-free properties of \(\bu \) and \(\tilde{\bu}\), it follows that
\begin{equation}\label{wsu:equ3}
\begin{aligned}
&\int_{Q_\tau} (\bu - \tilde{\bu}) \cdot \diver( \tilde{\bu}\otimes \tilde{\bu})   \\&= \int_{Q_\tau} \left[ (\bu - \tilde{\bu}) \cdot \tilde{\bu}\diver(\tilde{\bu}) + (\bu - \tilde{\bu} ) \otimes  \tilde{\bu} : \nabla \tilde{\bu} \right]
\\
&= \frac{1}{2}\int_{Q_\tau }(\bu - \tilde{\bu} ) \cdot  \nabla (|\tilde{\bu}|^2) 
\\&= 0.
\end{aligned}
\end{equation}
Substituting  (\ref{wsu:equ2}) and (\ref{wsu:equ3}) into (\ref{wsu:equ1}) yields
\begin{equation}\label{wsu:equ4}
\begin{aligned}
&\Big[ \int_\Omega \left( \frac{|\bu - \tilde{\bu}|^2}{2} + \theta\right) \Big]_{t=0}^{t = \tau} + \int_{Q_\tau}\left[   - \BS : \BD \tilde{\bu} - \tilde{\BS} : \BD \bu + \tilde{\BS} : \BD \tilde{\bu} \right] 
\\
&\leq
- \int_{Q_\tau} ( \bu - \tilde{\bu}) \otimes ( \bu - \tilde{\bu}) :\BD \tilde{\bu} .
\end{aligned}
\end{equation}

Next, we need to use of the entropy inequality for the weak solution and entropy balance for the classical solution. Testing in the entropy inequality (\ref{eq:weakPDE_entropy}) against \( \tilde{\theta}\) and using the \texblk{identification} \( S  = \log\theta\), we obtain
\begin{align*}
&\Big[ \int_\Omega \tilde{\theta} \texblk{S} \Big]_{t = 0}^{t = \tau} 
- \int_{Q_\tau} \left[ \texblk{S}\partial_t\tilde{\theta} 
+ \texblk{S}\bu \cdot \nabla \tilde{\theta} 
+ \frac{\bq}{\theta} \cdot \nabla \tilde{\theta}\right] 
\\
& \geq \texblk{-} \int_{Q_\tau} \tilde{\theta}\bq \cdot \frac{\nabla \theta}{\theta^2}  
+ \int_{Q_\tau} \frac{\BS : \BD \bu }{\theta} \tilde{\theta} .
\end{align*}
Adding this   to (\ref{wsu:equ4}), we get
\begin{equation}\label{wsu:equ5}
\begin{aligned}
&\Big[ \int_\Omega\left( \frac{|\bu - \tilde{\bu}|^2}{2} + \theta - \tilde{\theta}\log\theta\right) \Big]_{t = 0}^{ t = \tau } 
+ \int_{Q_\tau}\left[  - \frac{\tilde{\theta}}{\theta} \bq \cdot \frac{\nabla \theta}{\theta} + \frac{\tilde{\theta}}{\theta} \bq \cdot \frac{\nabla \tilde{\theta}}{\tilde{\theta}} \right]
\\
&\quad
+ \int_{Q_\tau} \left[ \frac{\tilde{\theta}}{\theta}\BS : \BD \bu - \BS : \BD \tilde{\bu} - \tilde{\BS} : \BD \bu  + \tilde{\BS} : \BD \tilde{\bu} \right] 
\\
&\leq - \int_{Q_\tau} \left[ ( \bu - \tilde{\bu}) \otimes ( \bu - \tilde{\bu} ) : \BD \tilde{\bu} +  \log\theta \partial_t\tilde{\theta} +  \log\theta \bu \cdot \nabla \tilde{\theta}\right].
\end{aligned}
\end{equation}
%\VP{I think that reviewer was incorrect in suggesting that the blue minus sign should be a minus sign in the above, simply comparing it to the previous line...?}
%\AG{The reviewer was correct, two signs were wrong, but I changed one more I shouldn't have}
We want to add terms to the left-hand side of \eqref{wsu:equ5} in order to have a term involving the relative energy.
Recalling the definition of the relative energy \( \mathcal{E}(\bu, \theta\,|\,\tilde{\bu}, \tilde{\theta}) \), we consider 
\begin{align*}
\Big[ \int_\Omega \left( - \tilde{\theta} +  \tilde{\theta}\log\tilde{\theta}\right)  \Big]_{t = 0}^{ t = \tau} = \int_{Q_\tau} \partial_t \tilde{\theta}\log\tilde{\theta} .
\end{align*}
Adding this into (\ref{wsu:equ5}), we deduce that
%\begin{equation}\label{wsu:equ9}
%\begin{aligned}
%&\Big[ \int_\Omega\mathcal{E}(\bu, \theta\,|\, \tilde{\bu}, \tilde{\theta}) \,\mathrm{d}x \Big]_{t = 0}^{t = \tau} + \int_{Q_\tau} \frac{\tilde{\theta}}{\theta} \BS : \BD \bu - \BS : \BD \tilde{\bu} - \tilde{\BS} : \BD \bu + \tilde{\BS} : \BD \tilde{\bu} \,\mathrm{d}x\,\mathrm{d}t
%\\
%&\quad
%- \int_{Q_\tau}  \frac{\tilde{\theta}}{\theta} \bq \cdot \frac{\nabla \theta}{\theta} - \frac{\tilde{\theta}}{\theta} \bq \cdot \frac{\nabla \tilde{\theta}}{\tilde{\theta}} \,\mathrm{d}x\,\mathrm{d}t
%\\
%&\leq - \int_{Q_\tau} (\bu - \tilde{\bu}) \otimes ( \bu - \tilde{\bu}) : \BD \tilde{\bu} + c_v \log
%\theta \partial_t \tilde{\theta} + c_v \log\theta \bu \cdot \nabla \tilde{\theta} - c_v \partial_t\tilde{\theta} \log\tilde{\theta} \,\mathrm{d}x\,\mathrm{d}t.
%\end{aligned}
%\end{equation}
\begin{equation}\label{wsu:equ9}
\begin{aligned}
&\Big[ \int_\Omega\mathcal{E}(\bu, \theta\,|\, \tilde{\bu}, \tilde{\theta})  \Big]_{t = 0}^{t = \tau} 
+ \int_{Q_\tau}\left[  \frac{\tilde{\theta}}{\theta} \BS \fp \BD \bu - \BS \fp \BD \tilde{\bu} - \tilde{\BS} \fp \BD \bu + \tilde{\BS} \fp \BD \tilde{\bu} \right]
\\
&\quad
- \int_{Q_\tau} \left[  \frac{\tilde{\theta}}{\theta} \bq \cdot \frac{\nabla \theta}{\theta} - \frac{\tilde{\theta}}{\theta} \bq \cdot \frac{\nabla \tilde{\theta}}{\tilde{\theta}} \right]
\\
&\leq - \int_{Q_\tau} \left[ (\bu - \tilde{\bu}) \otimes ( \bu - \tilde{\bu}) \fp \BD \tilde{\bu} +  \log
\theta \partial_t \tilde{\theta} +  \log\theta \bu \cdot \nabla \tilde{\theta} -  \partial_t\tilde{\theta} \log\tilde{\theta}\right]  .
\end{aligned}
\end{equation}
Using the energy balance (\ref{eq:PDE_energy}) for the classical solution, multiplying   by \( \log
\tilde{\theta} - \log\theta\) and integrating over \( Q_\tau\), we see that the second and fourth terms on the right-hand side of \eqref{wsu:equ9} can be rewritten as
\begin{equation}\label{wsu:equ6}
\int_{Q_\tau}  \partial_t\tilde{\theta} \big( \log
\tilde{\theta} - \log\theta\big) %\,\mathrm{d}x\,\mathrm{d}t
= \int_{Q_\tau} \big( \log\tilde{\theta} - \log\theta\big) \big[ \tilde{\BS} :\BD \tilde{\bu} - \diver (\tilde{\theta}\tilde{\bu}  + \tilde{\bq} ) \big].
\end{equation}
Noting that \( \log\theta\in L^2(0, T; W^{1,2}(\Omega) ) \) with \( \nabla \log
\theta \) identified by \( \theta^{-1}\nabla \theta\), we integrate  by parts and use the Neumann boundary condition on the terms involving the flux \( \tilde{\bq} \)  to deduce that
\begin{equation}\label{wsu:equ7}
 \int_{Q_\tau} \big( \log\tilde{\theta} - \log\theta\big) \diver\tilde{\bq}  = \int_{Q_\tau}\Big( \frac{\nabla \theta}{\theta} - \frac{\nabla \tilde{\theta}}{\tilde{\theta}}\Big) \cdot \tilde{\bq}  .
\end{equation}
As a result of  the incompressibility constraint on \( \tilde{\bu}\), we also see that
\begin{equation}\label{wsu:equ8}
\int_{Q_\tau} \big( \log\tilde{\theta} - \log\theta\big) \diver (\tilde{\theta}\tilde{\bu})   = \int_{Q_\tau}  \big( \log\tilde{\theta} - \log\theta\big) \tilde{\bu} \cdot \nabla \tilde{\theta}  .
\end{equation}
For the remaining term on the right-hand side of (\ref{wsu:equ6}), first write
\begin{align*}
\log\theta - \log\tilde{\theta} = \frac{\theta - \tilde{\theta}}{\tilde{\theta}} - \frac{(\theta - \tilde{\theta})^2}{\xi_{\theta, \tilde{\theta}}^2},
\end{align*}
where \( \xi_{\theta, \tilde{\theta}} \in [ \min\{ \theta, \, \tilde{\theta}\}, \max \{ \theta, \tilde{\theta}\}]\), by an application of Taylor's theorem. Recalling that \( \tilde{\theta}\) is uniformly bounded away from \(0 \) and is bounded above, it follows that
\begin{align*}
  \frac{(\theta - \tilde{\theta})^2}{\xi_{\theta, \tilde{\theta}}^2}  \leq C \big( |\log\theta| + \theta + 1\big) .
\end{align*}
Since \( \theta\) is bounded away from \(0 \),  
\begin{align*}
 \frac{(\theta - \tilde{\theta})^2}{\xi_{\theta, \tilde{\theta}}^2} \leq C |\theta - \tilde{\theta}|^2.
\end{align*}
Substituting the above, (\ref{wsu:equ7}) and (\ref{wsu:equ8}) in (\ref{wsu:equ6}), we use the result in the inequality (\ref{wsu:equ9}) and  deduce that

 \begin{equation}\label{wsu:equ10}
 \begin{aligned}
 &\Big[ \int_\Omega \mathcal{E}(\bu, \theta\,|\, \tilde{\bu}, \tilde{\theta}) \Big]_{t = 0}^{t = \tau} + \int_{Q_\tau}\Big[  \frac{\tilde{\theta}}{\theta}\BS \fp \BD \bu + \frac{\theta}{\tilde{\theta}}\tilde{\BS}\fp \tilde{\BD} - \tilde{\BS} \fp \BD \bu - \BS \fp \BD \tilde{\bu}\Big]
 \\
 &\quad
 - \int_{Q_\tau}\Big[  \frac{\tilde{\theta}}{\theta}\frac{\nabla \theta}{\theta}\cdot \bq
- \frac{\tilde{\theta}}{\theta}\frac{\nabla \tilde{\theta}}{\tilde{\theta}} \cdot \bq 
- \frac{\nabla \theta}{\theta}\cdot \tilde{\bq}  
+ \frac{\nabla \tilde{\theta}}{\tilde{\theta}}\cdot \tilde{\bq}\Big] 
\\
&\leq
- \int_{Q_\tau} \Big[ (\bu - \tilde{\bu}) \otimes ( \bu - \tilde{\bu}) \fp \BD \tilde{\bu} +  \log\theta\bu \cdot \nabla \tilde{\theta}  + (\log\tilde{\theta} - \log\theta
) \tilde{\bu}\cdot \nabla \tilde{\theta}\Big] 
\\
&\quad + \int_{Q_\tau} \frac{( \theta - \tilde{\theta})^2}{\xi_{\theta, \tilde{\theta}}^2} \tilde{\BS} \fp \BD \tilde{\bu}.
 \end{aligned}
 \end{equation}
  For the second and third terms on the right-hand side of (\ref{wsu:equ10}), we have
  \begin{align*}
 & - \int_{Q_\tau}\Big[  \log\theta\bu \cdot \nabla \tilde{\theta} + (\log\tilde{\theta} - \log\theta) \tilde{\bu}\cdot \nabla \tilde{\theta} \Big] 
  \\
  &= - \int_{Q_\tau} \Big[  \log\theta ( \bu - \tilde{\bu}) \cdot \nabla \tilde{\theta} +  \log\tilde{\theta}\tilde{\bu} \cdot \nabla \tilde{\theta} \Big] 
  \\
  &= - \int_{Q_\tau}  \log\theta ( \bu - \tilde{\bu}) \cdot \nabla \tilde{\theta} 
  \\
  &=  - \int_{Q_\tau}  \big( \log\theta - \log\tilde{\theta}\big) (\bu - \tilde{\bu}) \cdot \nabla \tilde{\theta}.
  \end{align*}
  Substituting this into (\ref{wsu:equ10}), we get
\begingroup
\allowdisplaybreaks
  \begin{equation}\nonumber
  \begin{aligned}
    &\Big[ \int_\Omega\mathcal{E}( \bu, \theta\,|\, \tilde{\bu}, \tilde{\theta})  \Big]_{t = 0}^{t = \tau} + \int_{Q_\tau}\Big[  \frac{\tilde{\theta
 }}{\theta} \BS : \BD \bu + \frac{\theta}{\tilde{\theta}} \tilde{\BS} : \BD\tilde{\bu } - \tilde{\BS}: \BD \bu - \BS : \BD \tilde{\bu} \Big] 
 \\
 &\quad - \int_{Q_\tau} \Big[ \frac{\tilde{\theta}}{\theta}\frac{\nabla \theta}{\theta}\cdot \bq
- \frac{\tilde{\theta}}{\theta}\frac{\nabla \tilde{\theta}}{\tilde{\theta}} \cdot \bq - \frac{\nabla \theta}{\theta}\cdot \tilde{\bq}  + \frac{\nabla \tilde{\theta}}{\tilde{\theta}}\cdot \tilde{\bq}\Big] 
\\
&\leq
- \int_{Q_\tau} ( \bu - \tilde{\bu }) \otimes (\bu - \tilde{\bu} ) : \BD \tilde{\bu}  
\\&\quad
- \int_{Q_\tau}  \big(\log\theta - \log\tilde{\theta}\big) (\bu - \tilde{\bu} ) \cdot \nabla \tilde{\theta} 
+ \int_{Q_\tau} \frac{(\theta - \tilde{\theta})^2}{\xi_{\theta, \tilde{\theta}}^2}  \tilde{\BS} : \BD \tilde{\bu} 
\\
&\leq \|\BD\tilde{\bu}\|_{L^\infty(Q_\tau)} \|\bu - \tilde{\bu}\|_{L^2(Q_\tau)}^2 + \|\tilde{\BS}\|_{L^\infty(Q_\tau)} \|\BD \tilde{\bu}\|_{L^\infty(Q_\tau)} \int_{Q_\tau} \frac{|\theta - \tilde{\theta}|^2 }{\xi_{\theta, \tilde{\theta}}^2} 
\\&\quad
+ \|\nabla \tilde{\theta}\|_{L^\infty(Q_\tau)} \|\log\theta - \log\tilde{\theta}\|_{L^2(Q_\tau)} \|\bu - \tilde{\bu}\|_{L^2(Q_\tau)}.
\end{aligned}
\end{equation}
\endgroup
Next, we notice that
\begin{align*}
\|\log\theta - \log\tilde{\theta}\|_{L^2(Q_\tau)}^2
&= \int_{Q_\tau} \left( \left[ |\log\theta - \log\tilde{\theta}|^2 \right]_{ess} + \left[ |\log\theta  - \log\tilde{\theta}|^2 \right]_{res} \right) 
\\
&\leq C \int_{Q_\tau}\left(  \left[ |\theta- \tilde{\theta}|^2 \right]_{ess} +  \left[ (\log\theta)^2 + (\log\tilde{\theta})^2 \right]_{res}\right) 
\\
&\leq C \int_{Q_\tau}\left(  \left[ |\theta- \tilde{\theta}|^2 \right]_{ess} + \left[ 1 + \theta\right]_{res}\right) 
\\
&\leq C \int_{Q_\tau} \mathcal{E}(\bu, \theta\,|\, \tilde{\bu}, \tilde{\theta}) ,
\end{align*}
where we use the fact that \( \theta\) is uniformly bounded away from \(0 \). Similarly, we have
\begin{align*}
\|\bu-\tilde{\bu}\|_{L^2(Q_\tau)}^2 \leq C  \int_{Q_\tau} \mathcal{E}(\bu, \theta\,|\, \tilde{\bu}, \tilde{\theta}) .
\end{align*}
Furthermore, we see that
\begin{align*}
\int_{Q_\tau} \frac{|\theta - \tilde{\theta}|^2 }{\xi_{\theta, \tilde{\theta}}^2} 
&\leq \int_{Q_\tau}\Big(  \Big[ \frac{|\theta - \tilde{\theta}|^2 }{\xi_{\theta, \tilde{\theta}}^2}\Big]_{ess} + \Big[ \frac{|\theta - \tilde{\theta}|^2 }{\xi_{\theta, \tilde{\theta}}^2}\Big]_{res}\Big) 
\\
&\leq C \int_{Q_\tau}\Big(  \big[ |\theta - \tilde{\theta}|^2\big]_{ess} + \big[ 1 + \theta + |\log\theta|\big]_{res} \Big) 
\\
&\leq C\int_{Q_\tau} \mathcal{E}(\bu, \theta\,|\,\tilde{\bu}, \tilde{\theta}) .
\end{align*}
It follows that
\begin{equation}\label{wsu:equ12}
\begin{aligned}
&\Big[ \int_\Omega\mathcal{E}( \bu, \theta\,|\, \tilde{\bu}, \tilde{\theta})   \Big]_{t = 0}^{t = \tau} + \int_{Q_\tau}\Big[  \frac{\tilde{\theta
 }}{\theta} \BS : \BD \bu + \frac{\theta}{\tilde{\theta}} \tilde{\BS} : \BD\tilde{\bu } - \tilde{\BS}: \BD \bu - \BS : \BD \tilde{\bu} \Big] 
 \\
 &\quad - \int_{Q_\tau}\Big[  \frac{\tilde{\theta}}{\theta}\frac{\nabla \theta}{\theta}\cdot \bq
- \frac{\tilde{\theta}}{\theta}\frac{\nabla \tilde{\theta}}{\tilde{\theta}} \cdot \bq - \frac{\nabla \theta}{\theta}\cdot \tilde{\bq}  + \frac{\nabla \tilde{\theta}}{\tilde{\theta}}\cdot \tilde{\bq}\Big] 
\\
&\leq C \int_{Q_\tau} \mathcal{E}( \bu, \theta\,|\, \tilde{\bu}, \tilde{\theta}) .
\end{aligned}
\end{equation}
We aim   to apply Gronwall's inequality to (\ref{wsu:equ12}) to deduce that
  \begin{equation}\label{wsu:equ13}
    \int_{\Omega}\mathcal{E}( \bu, \theta\,|\, \tilde{\bu}, \tilde{\theta})(t)  \leq \Big( \int_\Omega\mathcal{E}( \bu_0, \theta_0\,|\, \tilde{\bu}_0, \tilde{\theta}_0) \Big) \mathrm{e}^{Ct},
  \end{equation}
where \( C\) is a constant that is independent of \( t\). The required weak-strong uniqueness follows from this immediately by noticing that the right-hand side vanishes when the initial data coincide. However,  the second and third integrals on the left-hand side of (\ref{wsu:equ12}) are not necessarily non-negative. Thus we cannot apply Gronwall's inequality at present.

However, \texblk{under} appropriate assumptions on \( \tilde{\kappa}\) and  \( \Scr\) as stated in Assumption \ref{as:WSU_assumptions}, we are able to bound the integrals from below in a suitable way such that (\ref{wsu:equ13}) can be obtained.
Using Fourier's law, the constitutive relation concerning the heat flux term, we see that
\begin{align*}
&- \int_{Q_\tau} \Big[ \frac{\tilde{\theta}}{\theta}\frac{\nabla \theta}{\theta}\cdot \bq
- \frac{\tilde{\theta}}{\theta}\frac{\nabla \tilde{\theta}}{\tilde{\theta}} \cdot \bq - \frac{\nabla \theta}{\theta}\cdot \tilde{\bq}  + \frac{\nabla \tilde{\theta}}{\tilde{\theta}}\cdot \tilde{\bq}\Big] %\,\mathrm{d}x\,\mathrm{d}t
\\
&= \int_{Q_\tau}\Big[  \frac{\tilde{\theta}}{\theta} \frac{\nabla \theta}{\theta} \cdot \tilde{\kappa} ( \theta) \nabla \theta - \frac{\tilde{\theta}}{\theta} \frac{\nabla \tilde{\theta}}{\tilde{\theta}} \cdot \tilde{\kappa}( \theta) \nabla \theta - \frac{\nabla \theta}{\theta} \cdot \tilde{\kappa}(\tilde{\theta}) \nabla \tilde{\theta} + \frac{\nabla \tilde{\theta}}{\tilde{\theta}} \cdot \tilde{\kappa}(\tilde{\theta}) \nabla \tilde{\theta}\Big]%\,\mathrm{d}x\,\mathrm{d}t
\\
&= \int_{Q_\tau} \Big[ \tilde{\kappa}( \theta) \tilde{\theta} \Big| \frac{\nabla \theta}{\theta} - \frac{\nabla \tilde{\theta}}{\tilde{\theta}}\Big|^2 - \tilde{\kappa}( \theta) \frac{|\nabla \tilde{\theta}|^2 }{\tilde{\theta}} + 2\tilde{\kappa}( \theta) \frac{\nabla \theta}{\theta} \cdot \nabla \tilde{\theta} + \tilde{\kappa}( \tilde{\theta}) \frac{|\nabla \tilde{\theta}|^2}{\tilde{\theta}}  
\\&\quad
- \tilde{\kappa}( \theta) \frac{\nabla \theta}{\theta} \cdot \nabla \tilde{\theta}
- \tilde{\kappa}( \tilde{\theta}) \frac{\nabla \theta}{\theta} \cdot \nabla \tilde{\theta} \Big] %\,\mathrm{d}x\,\mathrm{d}t
\\
&= \int_{Q_\tau} \Big[ \tilde{\kappa}(\theta) \tilde{\theta}\Big| \frac{\nabla \theta}{\theta} - \frac{\nabla \tilde{\theta}}{\tilde{\theta}}\Big|^2 + \frac{|\nabla \tilde{\theta}|^2}{\tilde{\theta}} \big[ \tilde{\kappa}( \tilde{\theta}) - \tilde{\kappa}(\theta) \big] + \frac{\nabla \theta}{\theta}\cdot \nabla \tilde{\theta} \big[ \tilde{\kappa}(\theta) - \tilde{\kappa}( \tilde{\theta}) \big]\Big]  %\,\mathrm{d}x\,\mathrm{d}t
\\
&= \int_{Q_\tau} \Big[ \tilde{\kappa}(\theta) \tilde{\theta}\Big| \frac{\nabla \theta}{\theta} - \frac{\nabla \tilde{\theta}}{\tilde{\theta}}\Big|^2 + \nabla \tilde{\theta}\cdot \Big[ \frac{\nabla \tilde{\theta}}{\tilde{\theta}} - \frac{\nabla \theta}{\theta} \Big] \big[ \tilde{\kappa}( \tilde{\theta}) - \tilde{\kappa}( \theta) \big]\Big] % \,\mathrm{d}x\,\mathrm{d}t
\\
&\geq \int_{Q_\tau}\Big[ \frac{\tilde{\kappa}(\theta)\tilde{\theta}}{2} \Big| \frac{\nabla \theta}{\theta} - \frac{\nabla \tilde{\theta}}{\tilde{\theta}}\Big|^2 - \frac{|\nabla \tilde{\theta}|^2 }{2\tilde{\kappa}(\theta) \tilde{\theta}} |\tilde{\kappa}( \theta) - \tilde{\kappa}(\tilde{\theta}) |^2\Big]  %\,\mathrm{d}x\,\mathrm{d}t.
\end{align*}
Recalling that \( \tilde{\kappa} \) grows like \( \theta^{\frac{1}{2}}\) for large \( \theta\) and is locally Lipschitz, and using the boundedness of \( \tilde{\theta}\), we see that
\begin{align*}
\int_{Q_\tau}\frac{|\nabla \tilde{\theta}|^2 }{2\tilde{\kappa}(\theta) \tilde{\theta}} |\tilde{\kappa}( \theta) - \tilde{\kappa}(\tilde{\theta}) |^2% \,\mathrm{d}x\,\mathrm{d}t
&\leq C \int_{Q_\tau}\Big[  \big[ |\tilde{\kappa}(\theta) - \tilde{\kappa}(\tilde{\theta})|^2 \big]_{ess} + \big[ |\tilde{\kappa} ( \theta) - \tilde{\kappa}( \tilde{\theta}) |^2 \big]_{res}\Big]  %\,\mathrm{d}x\,\mathrm{d}t
\\
&\leq C\int_{Q_\tau} \Big[ \big[ |\theta - \tilde{\theta}|^2 \big]_{ess} + \big[ 1 + \theta \big]_{res}\Big]  %\,\mathrm{d}x\,\mathrm{d}t
\\
&\leq \int_{Q_\tau} \mathcal{E}( \bu, \theta\,|\, \tilde{\bu}, \tilde{\theta}). %\,\mathrm{d}x\,\mathrm{d}t.
\end{align*}
Substituting this bound into (\ref{wsu:equ12}), we deduce that
\begin{equation}\label{wsu:equ14}
\begin{aligned}
&\Big[ \int_\Omega\mathcal{E}(\bu, \theta\,|\, \tilde{\bu}, \tilde{\theta}) \Big]_{t = 0}^{t = \tau} + \int_{Q_\tau} \Big[  \frac{\tilde{\theta
 }}{\theta} \BS : \BD \bu + \frac{\theta}{\tilde{\theta}} \tilde{\BS} : \BD\tilde{\bu } - \tilde{\BS}: \BD \bu - \BS : \BD \tilde{\bu} \Big] 
 \\
 &\quad
% + \int_{Q_\tau} \frac{|\nabla \tilde{\theta}|^2}{2\tilde{\kappa}(\theta) \tilde{\theta}} |\tilde{\kappa}(\theta) - \tilde{\kappa}( \tilde{\theta}) |^2 \,\mathrm{d}x\,\mathrm{d}t
+\texblk{\int_{Q_\tau}\frac{\tilde{\kappa}(\theta)\tilde{\theta}}{2} \Big| \frac{\nabla \theta}{\theta} - \frac{\nabla \tilde{\theta}}{\tilde{\theta}}\Big|^2}
 \\
 &\leq C\int_{Q_\tau} \mathcal{E}(\bu, \theta\,|\, \tilde{\bu}, \tilde{\theta}).
\end{aligned}
\end{equation}

To deal with the second integral on the left-hand side of \eqref{wsu:equ14} we will now make use of the strong monotonicity condition \eqref{eq:strong_monotonicity}. Considering the integrand in the aforementioned term, we write
\begin{align*}
&\frac{\tilde{\theta}}{\theta}\BS : \BD\bu + \frac{\theta}{\tilde{\theta}} \tilde{\BS}: \BD \tilde{\bu} - \tilde{\BS} : \BD \bu - \BS : \BD \tilde{\bu}
\\
&
= \frac{\tilde{\theta}}{\theta}( \BS - \tilde{\BS}) : (\BD \bu - \BD \tilde{\bu} ) + \tilde{\theta} \Big( \frac{1}{\theta} - \frac{1}{\tilde{\theta}} \Big) \tilde{\BS} : \BD \bu
+ \tilde{\theta} \Big( \frac{1}{\theta} - \frac{1}{\tilde{\theta}}\Big) \BS : \BD \tilde{\bu} 
\\&\quad
+ \theta\Big( \frac{1}{\tilde{\theta}} - \frac{1}{\theta}\Big) \tilde{\BS} : \BD \tilde{\bu} 
+ \tilde{\theta} \Big( \frac{1}{\tilde{\theta}} - \frac{1}{\theta}\Big) \tilde{\BS} : \BD \tilde{\bu}
\\
& 
= \frac{\tilde{\theta}}{\theta} ( \BS - \Scr(\BD\tilde{\bu},\theta)) : (\BD \bu - \BD \tilde{\bu})
+ \frac{\tilde{\theta}}{\theta} ( \Scr(\BD\tilde{\bu},\theta) - \tilde{\BS}) : (\BD \bu - \BD \tilde{\bu})
\\
&\quad
+ \tilde{\theta} \Big( \frac{1}{\theta} - \frac{1}{\tilde{\theta}} \Big) \tilde{\BS} : (\BD \bu - \BD \tilde{\bu} )
+ \tilde{\theta}\Big( \frac{1}{\theta} - \frac{1}{\tilde{\theta}} \Big) \BD \tilde{\bu} : (\BS - \tilde{\BS})
\\
&\quad + \Big( \frac{1}{\theta} - \frac{1}{\tilde{\theta}}\Big) (\tilde{\theta} - \theta) \tilde{\BS} :\BD \tilde{\bu}
\\
&\geq \frac{\tilde{\theta}}{\theta} c \big( |\BS - \tilde{\BS}|^2 + |\BD \bu - \BD \tilde{\bu}|^2 \big)
+ \frac{\tilde{\theta}}{\theta} ( \Scr(\BD\tilde{\bu},\theta) - \Scr(\BD\tilde{\bu},\tilde{\theta})) : (\BD \bu - \BD \tilde{\bu})
\\
&\quad
+ \tilde{\theta} \Big( \frac{1}{\theta} - \frac{1}{\tilde{\theta}}\Big) \tilde{\BS} : (\BD \bu - \BD \tilde{\bu})
+ \tilde{\theta}\Big( \frac{1}{\theta} - \frac{1}{\tilde{\theta}} \Big) \BD \tilde{\bu} : (\BS - \tilde{\BS})
\\ 
&\quad +\Big( \frac{1}{\theta} - \frac{1}{\tilde{\theta}}\Big) (\tilde{\theta} - \theta) \tilde{\BS}:\BD\tilde{\bu}.
\end{align*}
The second term in the right-hand side can be dealt with by means of Young's inequality, which yields
\begin{align*}
\frac{\tilde{\theta}}{\theta} ( \Scr(\BD\tilde{\bu},\theta) - \Scr(\BD\tilde{\bu},\tilde{\theta})) \fp (\BD \bu - \BD \tilde{\bu})
&\leq \epsilon \frac{\tilde{\theta}}{\theta} |\BD\bu - \BD\tilde{\bu}|^2
\\
&\quad
+ C(\epsilon)\frac{\tilde{\theta}}{\theta}|\Scr(\BD\tilde{\bu},\theta) - \Scr(\BD\tilde{\bu},\tilde{\theta})|^2,
\end{align*}
where \( \epsilon >0 \) is a constant to be fixed sufficiently small later on.
In a similar manner, for the third and fourth terms on the right-hand side, we have
\begin{align*}
&\tilde{\theta}\Big( \frac{1}{\theta} - \frac{1}{\tilde{\theta}}\Big) \tilde{\BS} : (\BD \bu - \BD \tilde{\bu}) + \tilde{\theta} \Big( \frac{1}{\theta} - \frac{1}{\tilde{\theta}}\Big) \BD \tilde{\bu} : (\BS - \tilde{\BS})
\\
&
= \sqrt{\tilde{\theta}\theta}\Big( \frac{1}{\theta} - \frac{1}{\tilde{\theta}}\Big) \sqrt{\frac{\tilde{\theta}}{\theta}}
\tilde{\BS}: (\BD \bu - \BD \tilde{\bu} ) + \sqrt{\tilde{\theta}\theta}\Big( \frac{1}{\theta} - \frac{1}{\tilde{\theta}}\Big) \sqrt{\frac{\tilde{\theta}}{\theta}}
\BD \tilde{\bu}: (\BS - \tilde{\BS})
\\
&\leq \epsilon \frac{\tilde{\theta}}{\theta} \Big( |\BD \bu - \BD \tilde{\bu }|^2  + |\BS - \tilde{\BS}|^2  \Big)
+ C(\epsilon)\Big( \|\tilde{\BS}\|_{L^\infty(Q_\tau)}  + \|\BD \tilde{\bu}\|_{L^\infty(Q_\tau)}\Big)  \frac{|\theta - \tilde{\theta}|^2}{\theta\tilde{\theta}} .
\end{align*}
Choosing \( \epsilon = \frac{c}{4}\) in the two previous estimates, where \( c\) is the constant from the uniform monotonicity assumption, we deduce that
\begin{align*}
&\int_{Q_\tau} \Big[ \frac{\tilde{\theta}}{\theta}\BS : \BD\bu + \frac{\theta}{\tilde{\theta}} \tilde{\BS}: \BD \tilde{\bu} - \tilde{\BS} : \BD \bu - \BS : \BD \tilde{\bu}  \Big] 
\\
&\geq \int_{Q_\tau} \frac{c}{2}\frac{\tilde{\theta}}{\theta} \Big( |\BS - \tilde{\BS}|^2 + |\BD \bu - \BD \tilde{\bu}|^2 \Big)  
+ \int_{Q_\tau}\Big( \frac{1}{\theta} - \frac{1}{\tilde{\theta}}\Big) (\tilde{\theta} - \theta) \tilde{\BS}: \BD \tilde{\bu} 
\\&\quad
- C \int_{Q_\tau}\frac{|\theta - \tilde{\theta}|^2}{\theta\tilde{\theta}}  
- C \int_{Q_\tau} \frac{\tilde{\theta}}{\theta}|\Scr(\BD\tilde{\bu},\theta) - \Scr(\BD\tilde{\bu},\tilde{\theta})|^2 .
\end{align*}
\texblk{The second term on the right-hand side can be dropped since it is non-negative, thanks to the decreasing character of the function $x\mapsto \frac{1}{x}$ for $x>0$. This yields}
%Using the uniform boundedness of \( \tilde{\BS}\) and \( \BD \tilde{\bu}\), the second term on the right-hand side of the above can be absorbed into the third to yield
\begin{equation}\label{wsu:equ15}
    \begin{aligned}
        &\int_{Q_\tau}\Big[  \frac{\tilde{\theta}}{\theta}\BS : \BD\bu + \frac{\theta}{\tilde{\theta}} \tilde{\BS}: \BD \tilde{\bu} - \tilde{\BS} : \BD \bu - \BS : \BD \tilde{\bu} \Big] 
    \\
    &\geq \int_{Q_\tau} \frac{c}{2}\frac{\tilde{\theta}}{\theta} \Big( |\BS - \tilde{\BS}|^2 + |\BD \bu - \BD \tilde{\bu}|^2 \Big)  
    - C \int_{Q_\tau}\frac{|\theta - \tilde{\theta}|^2}{\theta\tilde{\theta}}  
    \\
    &\quad 
    - C \int_{Q_\tau} \frac{\tilde{\theta}}{\theta}|\Scr(\BD\tilde{\bu},\theta) - \Scr(\BD\tilde{\bu},\tilde{\theta})|^2 .
    \end{aligned}
\end{equation}
The first integral on the right-hand side is of the desired form. For the second integral, we have
\begin{align*}
\int_{Q_\tau} \frac{|\theta - \tilde{\theta}|^2 }{\theta\tilde{\theta}} 
&= \int_{Q_\tau} \Big( \Big[ \frac{|\theta - \tilde{\theta}|^2 }{\theta\tilde{\theta}} \Big]_{ess} + \Big[  \frac{|\theta - \tilde{\theta}|^2 }{\theta\tilde{\theta}} \Big]_{res}\Big) 
\\
&\leq C \int_{Q_\tau} \Big( [|\theta - \tilde{\theta}|^2 ]_{ess}+ [1 + \theta]_{res} \Big) 
\\
&\leq C \int_{Q_\tau} \mathcal{E}(\bu, \theta \,|\, \tilde{\bu}, \tilde{\theta}) ,
\end{align*}
making use of the fact that \( \theta\) and \( \tilde{\theta}\) are uniformly bounded away from 0. The final integral on the right-hand side of \eqref{wsu:equ15} can be estimated likewise, thanks to the Lipschitz condition \eqref{eq:loc_Lipschitz_S} and the growth assumption \eqref{eq:growth}. We argue in the following way: 
\begin{align*}
  &\int_{Q_\tau} \frac{\tilde{\theta}}{\theta}|\Scr(\BD\tilde{\bu},\theta) - \Scr(\BD\tilde{\bu},\tilde{\theta})|^2 
\\&\leq C \int_{Q_\tau}\Big(  \Big[ |\Scr(\BD\tilde{\bu},\theta) - \Scr(\BD\tilde{\bu},\tilde{\theta})|^2\Big]_{ess} +  \Big[ |\Scr(\BD\tilde{\bu},\theta) - \Scr(\BD\tilde{\bu},\tilde{\theta})|^2\Big]_{res}\Big) 
\\ &
\leq C(\|\BD\tilde{\bu}\|_{L^\infty(Q_\tau)})\int_{Q_\tau}\Big(  [|\theta - \tilde{\theta}|^2]_{ess} +  [1]_{ess} \Big) 
\\&
\leq C\int_{Q_\tau} \mathcal{E}(\bu, \theta\,|\, \tilde{\bu}, \tilde{\theta}) ,
\end{align*}
using the Lipschitz property of \( \Scr\) to deal with the essential part and growth condition on \( \Scr\) to bound the residual part. 

It follows that
\begin{align*}
&\Big[ \int_\Omega \mathcal{E}(\bu, \theta\,|\, \tilde{\bu}, \tilde{\theta})  \Big]_{t = 0}^{t = \tau} + \int_{Q_\tau} \frac{\tilde{\theta}}{\theta} \Big( |\BS - \tilde{\BS}|^2 + |\BD \bu - \BD \tilde{\bu}|^2 \Big)  
\\
%&\quad
%+ \int_{Q_\tau} \Big( \frac{1}{\theta} - \frac{1}{\tilde{\theta}}\Big) (\tilde{\theta} - \theta) \tilde{\BS}: \BD \tilde{\bu} \,\mathrm{d}x\,\mathrm{d}t
%\\
&\leq C\int_{Q_\tau} \mathcal{E}(\bu, \theta\,|\, \tilde{\bu}, \tilde{\theta}) .
\end{align*}
Applying Gronwall's inequality, we conclude the required result. We note that the existence of strong solutions for small data has been obtained in \cite{Amann1995} (see also \cite{Benes2011}); under such assumptions we would then have that the numerical approximations devised in Proposition \ref{exist:prop1} would converge to the strong solution.

Some examples of constitutive relations satisfying the assumptions required in this section, in particular the strong monotonicity \eqref{eq:strong_monotonicity}, include the Carreau--Yasuda constitutive relation for $r<2$, given by
\begin{equation}\label{eq:CarreauYasuda}
  \BS = \alpha(\theta) \BD\bu + \beta(\theta)(1+ \Gamma(\theta)|\BD\bu|^2)^{\frac{r-2}{2}}\BD\bu,
\end{equation}
where $\alpha$, $\beta$ and $\Gamma$ are locally Lipschitz continuous functions satisfying $0<c_1 \leq \alpha,\beta,\Gamma \leq c_2$, for two positive constants $c_1,c_2>0$. Although the Herschel--Bulkley constitutive relation \eqref{eq:HerschelBulkley} is not strongly monotone, in practice it is common to regularise it when performing numerical approximations, in order to deal with its non-differentiable character. We mention for example, that the following regularisation (introduced in \cite{Bulicek2020})
\begin{equation}
\BG(\BS - \varepsilon\BD\bu,\BD\bu-\varepsilon\BS) = \bm{0},
\end{equation}
which can be applied to the expression \eqref{eq:HB_implicit}, leads to a relation satisfying the strong monotonicity condition \eqref{eq:strong_monotonicity}, even if $r\neq 2$ \cite[Eq.\ 4.26]{Bulicek2020}.

\begin{rmk}
  The weak-strong uniqueness result for the isothermal system obtained in \cite{Abbatiello2019} only requires that $r>1$ and surprisingly does not impose the constitutive relation pointwise at the level of the dissipative weak solutions. One consequence is that the result has to account for possible concentrations in the convective term $\diver(\bu\otimes \bu)$. An extension of the results presented in this paper to a setting with such relaxed assumptions will be the subject of future work. In the current non-isothermal setting, the assumption $r>1$ would also lead to potential concentrations in the advective term for the entropy  $\diver(S\bu)$. On the other hand, since the constitutive relation would not need to be identified pointwise, a discrete parabolic Lipschitz truncation would not be necessary, which means that the penalty term (and so the index $k$) could be dropped from the approximation scheme.
\end{rmk}

\begin{rmk}
  The notion of dissipative weak solution introduced in Definition \ref{def:dissipative_weak_sols} is suitable for energetically isolated systems with $\bu |_{\partial\Omega}=\bm{0}$ and $\bm{q}|_{\partial\Omega}=\bm{0}$. The case of Dirichlet boundary conditions for the temperature $\theta|_{\partial\Omega}= \theta_b$ is therefore excluded; this problem has only recently been solved in the compressible case in \cite{Chaudhuri2021}. In this case it was necessary to modify the balance of total energy \eqref{eq:weakPDE_energy} and employ instead the so-called ballistic free energy. It seems plausible that these arguments carry over to the incompressible setting and will be the subject of future research.
\end{rmk}

\section*{Acknowledgements}
The authors would like to thank the anonymous referees, whose detailed comments and suggestions helped to increase the quality of the manuscript.

\bibliographystyle{elsarticle-num} 
\bibliography{bibliography}

\end{document}